\documentclass[11pt]{amsart}

\usepackage{amsmath,amssymb,amsthm,mathtools}
\usepackage{enumitem}
\usepackage{float}
\usepackage{pgfplots}
\usepackage[margin=1.1in]{geometry}
\usepackage[colorlinks=true,linkcolor=blue,citecolor=blue,urlcolor=blue]{hyperref}
\hypersetup{
  pdftitle={Modulating flow polynomials and parity bias in homological walk counting},
  pdfauthor={Ye Luo}
}
\pgfplotsset{compat=1.15}
\newtheorem{theorem}{Theorem}[section]
\newtheorem{proposition}[theorem]{Proposition}
\newtheorem{lemma}[theorem]{Lemma}
\newtheorem{corollary}[theorem]{Corollary}

\theoremstyle{definition}
\newtheorem{definition}[theorem]{Definition}
\newtheorem{example}[theorem]{Example}
\theoremstyle{remark}
\newtheorem{remark}[theorem]{Remark}

\DeclareMathOperator{\Tr}{Tr}
\DeclareMathOperator{\Spec}{Spec}

\newcommand{\CG}{C_G}
\newcommand{\DG}{D_G}
\newcommand{\Deta}{D_{G,\eta}}
\newcommand{\GF}{G[F]}

\newcommand{\NZ}{\operatorname{NZ}}
\newcommand{\pc}{\operatorname{pc}}
\newcommand{\wtN}{\widetilde N}

\newcommand{\bbR}{\mathbb R}
\newcommand{\bbZ}{\mathbb Z}
\newcommand{\bbC}{\mathbb C}
\newcommand{\ee}{\mathrm e}

\title[Modulating flow polynomials]{Modulating flow polynomials and parity bias in homological walk counting}

\author{Ye Luo}
\address{School of Informatics, Xiamen University, Xiamen, Fujian 361005, China}
\email{luoye@xmu.edu.cn}
\date{}
\makeatletter
\@ifundefined{subjclassname@2020}{%
  \@namedef{subjclassname@2020}{\textup{2020} Mathematics Subject Classification}%
}{}
\makeatother

\begin{document}

\begin{abstract}
Long weighted walks seem to sample mod-\(t\) flows uniformly, but at even moduli a
persistent parity bias survives.  Writing \(C_G(t)\) for the flow polynomial of a
connected graph \(G\) of genus \(g\), the proportion of nowhere-zero flows under
uniform sampling is \(C_G(t)/t^g\).  For arbitrary positive directed weights and
all moduli, we determine the exact leading corrections for both closed walks and
nonbacktracking closed walks, which we call \emph{circuits}.  Only the circuit
result requires \(g\ge2\).  For closed walks, the canonical character governs the
leading correction through the canonical modulating flow polynomial \(D_G\); for
circuits, the period character, when it exists, governs the leading correction
through the period modulating flow polynomial \(D_{G,\eta}\).  At
even moduli these give parity-dependent corrections: the two mechanisms coincide
on nonbipartite graphs and separate only in the bipartite circuit problem.  More
generally, each integral-homology character \(\zeta\) determines a
balanced-dichromatic specialization \(D_{G,\zeta}\).  These polynomials form a
character-indexed family whose finite-order members give Fourier coefficients of
the nowhere-zero flow set and whose retained edge sets are the central sets of a
gain semimatroid.  We further show that \(D_G\) is a cycle-matroid invariant and a
signed-Tutte specialization satisfying signed deletion-contraction.  Finally, an
integral nowhere-zero approximation theorem implies that asymptotically almost
every weighted circuit traverses every edge.
\end{abstract}

\maketitle

\section{Introduction}

Throughout, a graph means a finite undirected multigraph; multiple edges and loops
are allowed unless otherwise stated.  In this paper, a graph is \emph{bipartite}
if it is loopless and its vertex set can be partitioned into two classes such that
every edge joins the two classes.  In particular, every graph containing a loop is
nonbipartite.  Consider a graph \(G\) with vertex set
\[
   V(G)=\{v_1,\ldots,v_n\}
\]
and edge set
\[
   E(G)=\{e_1,\ldots,e_m\}.
\]
In the introduction, when the graph is fixed, we also write \(G=(V,E)\).  For
\(F\subseteq E\), let
\[
   \GF=(V,F)
\]
be the spanning subgraph with edge set \(F\), and let \(c(F)\) be the number of
connected components of \(\GF\).  We write
\[
   \beta(F)=|F|-|V|+c(F).
\]
This is the cycle rank, or first Betti number, of the spanning subgraph \(\GF\);
equivalently \(\beta(F)=\operatorname{rank} H_1(\GF;\bbZ)\), the number of
independent cycles in \(\GF\).  Thus a \(t\)-torsion cycle space supported on \(F\)
has size \(t^{\beta(F)}\).  We also write
\[
   \beta(G)=\beta(E)=|E|-|V|+c(E).
\]
When \(G\) is connected, \(c(E)=1\), and the genus of \(G\) is precisely this
Betti number:
\[
   g=m-n+1=|E|-|V|+1=\beta(G).
\]

The flow polynomial of \(G\), introduced by Tutte in his work on graph
polynomials~\cite{Tutte1954}, is a central object in the theory of nowhere-zero
flows; see, for example, the surveys and related developments
in~\cite{Shahmohamad2013,Dong2020FlowZeros,Onn2004,BreuerSanyal2012,Kochol2022Chains}.
It may be written in the form
\begin{equation}\label{eq:flow-poly}
   \begin{aligned}
   \CG(t)
      &=
      \sum_{F\subseteq E}(-1)^{|E|-|F|}
          t^{|F|-|V|+c(F)}  \\
      &=
      \sum_{F\subseteq E}(-1)^{|E|-|F|}
          t^{\beta(F)}.
   \end{aligned}
\end{equation}
For positive integers \(t\), this polynomial counts nowhere-zero flows over any
abelian group of order \(t\).  Modular, integral, complementary, and signed-graph
variants, as well as their orientation and reciprocity structures, have been studied
in~\cite{Chen2007Orientations,ChenYang2009,BreuerSanyal2012,Chen2013DualComplementary,RenQian2019}.
We recall the standard bridge and degree properties: \(\CG\equiv0\) if and only if
\(G\) has a bridge, whereas for bridgeless \(G\), the polynomial \(\CG\) is monic of
degree \(\beta(G)\)~\cite[Properties~1.5 and~1.7]{Shahmohamad2013}.  In particular,
if \(G\) is connected and bridgeless, then
\(\CG(t)/t^g=1+O_G(t^{-1})\) as \(t\to\infty\).

Walk counting has a different classical language.  The total number of closed walks
is a trace of an adjacency matrix, while reduced or nonbacktracking closed walks,
which we call \emph{circuits}, are encoded by the Hashimoto matrix
and by Ihara-type zeta functions; see, for example, Terras~\cite{Terras2011}.
Homological refinements of such trace formulas weight a
walk by its abelianization in the cycle space, as in homological spectral graph
theory~\cite{LuoRoyHSGT} and in the Eulerian counting paper~\cite{LuoEulerian}.
Thus flow polynomials and walk-counting zeta functions approach the same graph from
opposite directions: the former count static nowhere-zero circulations, while the
latter count dynamically generated circulations coming from long closed walks.

This raises a natural comparison problem.  If \(G\) is connected, then
\(H_1(G,\bbZ/t\bbZ)\) has size \(t^g\), so \(\CG(t)/t^g\) is the nowhere-zero ratio in
the uniform mod-\(t\) flow space.  A long closed walk, or a long circuit, also
determines a mod-\(t\) circulation by abelianization.  One might
therefore expect the same nowhere-zero ratio in this dynamic model.  In other words,
the most naive prediction is that long walks sample the mod-\(t\) flow space
uniformly enough that the ordinary flow polynomial is the whole leading answer.

The main point of this paper is that this expectation is almost, but not exactly,
correct, and that the discrepancy is itself governed by flow polynomials.  The
ordinary flow polynomial is the trivial Fourier coefficient in the relevant trace
formula.  For closed walks there is at most one further coefficient at spectral
radius: when \(G\) is nonbipartite and \(t\) is even, the canonical coefficient
survives.  It is obtained from Tutte's expansion by retaining precisely
the spanning subgraphs which contain no odd cycle.

\begin{samepage}
\begin{definition}[Canonical modulating flow polynomial]
\label{def:modulating-flow-poly}
The \emph{modulating flow polynomial associated with the canonical character}, or
simply the \emph{modulating flow polynomial}, of \(G\) is
\begin{equation}\label{eq:modulating-flow-poly}
   \begin{aligned}
   \DG(t)
      &=
      \sum_{\substack{F\subseteq E\\ \GF\ \mathrm{bipartite}}}
          (-1)^{|E|-|F|}
          t^{|F|-|V|+c(F)}  \\
      &=
      \sum_{\substack{F\subseteq E\\ \GF\ \mathrm{bipartite}}}
          (-1)^{|E|-|F|}
          t^{\beta(F)}.
   \end{aligned}
\end{equation}
\end{definition}
\end{samepage}

Thus \(\DG(t)\) is not an auxiliary error term.  It is a flow-type invariant with the
same cycle-rank weights as \(\CG(t)\), but with the bipartite part of the subgraph
lattice singled out.  This restriction is exactly what is seen by the
\emph{canonical character} \(\theta\); see \eqref{eq:disconnected-canonical-character} and
\cite[Lemma~3.1 and Definition~3.2]{LuoRoyHSGT}.  Below we prove that, for even
\(t\), \(\DG(t)\) is the signed count of nowhere-zero mod-\(t\) flows according to
its descent \(\theta_t\).

For circuits the canonical character does not by itself describe
every graph.  The edge dynamics has a period \(\nu_G\), and the new spectral
classification in homological spectral graph theory~\cite[Theorem~3.6]{LuoRoyHSGT}
identifies at most one additional extremal character.  The resulting
\emph{period character} \(\eta\) is defined for connected graphs in
\cite[Definition~3.5]{LuoRoyHSGT}; for disconnected graphs we use the componentwise
extension in \eqref{eq:componentwise-period-character}.  When it exists, it records
the parity of circuit length after division by \(\nu_G\).  It equals \(\theta\) on a
nonbipartite graph,
but on a bipartite graph it may fail to exist or may be a distinct nontrivial
character.  The corresponding Fourier coefficient leads to a second polynomial,
\(\Deta(t)\), which completes the nonbacktracking theory.
The period character and its spectral rotation are supplied by homological spectral
graph theory; the new step here is to identify the surviving nowhere-zero Fourier
coefficient as a flow-type polynomial and place it in the signed-Tutte framework.

The main contribution of the paper is the complete weighted asymptotic theorem
below.  It shows that each transfer operator has at most one nontrivial leading
Fourier coefficient, but that the distinguished coefficient is operator-dependent:
\(\theta\) governs vertex walks, whereas \(\eta\) governs circuits.
The two mechanisms coincide on nonbipartite graphs and separate only in the
bipartite edge problem.

To prepare for this result, we first investigate the two distinguished modulating
flow polynomials \(D_G\) and \(D_{G,\eta}\) through their Fourier interpretation and
parity properties.  They are the canonical and period members of a natural
character-indexed family.  For every full complex
character \(\zeta:H_1(G,\bbZ)\to\bbC^\ast\), retaining those spanning subgraphs on
whose integral homology \(\zeta\) is trivial gives a specialization of Zaslavsky's
balanced dichromatic polynomial~\cite{Zaslavsky1995BiasedIII}.
Section~\ref{sec:fourier-coefficients} denotes it
by \(D_{G,\zeta}(q)\) and establishes its Fourier interpretation.  If \(\zeta\) has
finite order \(N\), then, precisely at moduli divisible by \(N\), its evaluation is
the corresponding Fourier coefficient of the nowhere-zero flow set.  Homological
spectral graph theory then selects only the trivial and canonical characters for
vertex walks, and, in genus at least two, only the trivial and possible period
characters for circuits.  Genus one is the exceptional edge case,
where every unitary character has maximal spectral
radius~\cite[Theorem~3.6]{LuoRoyHSGT}.

The word ``modulating'' is meant literally.  For closed walks, the exceptional
leading ratio is
\[
   \frac{\CG(t)+(-1)^\ell\DG(t)}{t^g}.
\]
For circuits, whenever \(\eta\) exists and \(t\) is even, the
corresponding ratio along lengths \(\ell\) divisible by \(\nu_G\) is
\[
   \frac{\CG(t)+(-1)^{\ell/\nu_G}\Deta(t)}{t^g}.
\]
Although the total weighted mass and spectral radius vary with an arbitrary positive
directed weight, these coefficients do not.  They measure the leading discrepancy
between uniform sampling from the mod-\(t\) flow space and dynamic sampling by long
weighted walks or circuits.

We now describe the counting problem precisely.  Assigning a direction to each edge
gives an orientation \(\mathfrak o\) of \(G\).  For each edge
\(e_i\in E(G)\), let \(\mathbf e_i\) and \(\mathbf e_i^{-1}\) be the two oriented
edges, positively and negatively oriented with respect to \(\mathfrak o\).  Put
\[
   \mathbf E_{\mathfrak o}(G)=\{\mathbf e_1,\ldots,\mathbf e_m\}
\]
and
\[
   \mathbf E(G)
      =
      \{\mathbf e_1,\ldots,\mathbf e_m,
        \mathbf e_1^{-1},\ldots,\mathbf e_m^{-1}\}.
\]
Thus \(G_{\mathfrak o}=(V(G),\mathbf E_{\mathfrak o}(G))\) is the digraph determined
by the orientation \(\mathfrak o\), and \(G\) is its underlying graph.  For
\(\mathbf e\in\mathbf E(G)\), write \(\mathbf e(0)\) and \(\mathbf e(1)\) for its
initial and terminal vertices, and \(\mathbf e^{-1}\) for its reverse.  A loop
contributes two formal oriented loops, inverse to each other.  For a subgraph
\(H\subseteq G\), write \(\mathbf E_{\mathfrak o}(H)\) for the positively oriented
edges of \(H\), and write
\[
   E^c(H)=E(G)\setminus E(H),
   \qquad
   \mathbf E_{\mathfrak o}^c(H)
      =
      \mathbf E_{\mathfrak o}(G)\setminus \mathbf E_{\mathfrak o}(H).
\]
A \emph{walk} of length \(\ell\) is a sequence
\(P=\mathbf e_{i_1}\cdots\mathbf e_{i_\ell}\) of oriented edges such that
\(\mathbf e_{i_j}(1)=\mathbf e_{i_{j+1}}(0)\) for \(1\le j<\ell\).  Its length
is denoted by \(|P|=\ell\).  It is \emph{closed} if also
\(\mathbf e_{i_\ell}(1)=\mathbf e_{i_1}(0)\).  The walk is
\emph{nonbacktracking} if
\(\mathbf e_{i_{j+1}}\ne\mathbf e_{i_j}^{-1}\) for \(1\le j<\ell\); for a
closed walk we additionally impose the cyclic condition
\(\mathbf e_{i_1}\ne\mathbf e_{i_\ell}^{-1}\).  Such a nonbacktracking closed
walk is called a \emph{circuit}.  Throughout this paper, an unqualified
\emph{circuit} always has this dynamic meaning.  In the matroid-theoretic sense, a
inclusion-minimal dependent set is called a \emph{matroid circuit} (or a circuit of
the named matroid), and its graph-theoretic counterpart in a graphic matroid is
called a \emph{cycle}.  Throughout the dynamic counting statements, a circuit has
a distinguished initial oriented edge; it
is therefore rooted and oriented, and cyclic shifts and reversals are counted
separately, as in
\(\operatorname{Tr}(B_W^{\ell})\).

For a commutative ring \(R\), \(C_1(G,R)\) is the free \(R\)-module on
\(\mathbf E_{\mathfrak o}(G)\), \(C_0(G,R)\) is the free \(R\)-module on \(V(G)\),
and
\[
   \partial \mathbf e=\mathbf e(1)-\mathbf e(0)
      \qquad(\mathbf e\in\mathbf E_{\mathfrak o}(G)).
\]
We regard an oriented edge \(\mathbf e\in\mathbf E(G)\) as the corresponding signed
generator of \(C_1(G,R)\).  Since a graph has no two-cells,
\(H_1(G,R)=\ker\partial\).  The coboundary of a potential \(\varphi:V(G)\to R\) is
\[
   d\varphi(\mathbf e)=\varphi(\mathbf e(1))-\varphi(\mathbf e(0))
      \qquad(\mathbf e\in\mathbf E_{\mathfrak o}(G)).
\]
Finally, \(\mathbf 1_{\mathcal A}\) denotes the indicator of a condition or subset
\(\mathcal A\).

If
\[
   P=\mathbf e_{i_1}\mathbf e_{i_2}\cdots \mathbf e_{i_\ell}
\]
is a closed walk of length \(\ell\), with each \(\mathbf e_{i_j}\in\mathbf E(G)\),
its abelianization is the integral chain
\[
   P^{\mathrm{ab}}
      =
      \mathbf e_{i_1}+\mathbf e_{i_2}+\cdots+\mathbf e_{i_\ell}
      \in C_1(G,\bbZ).
\]
Since \(P\) is closed, \(\partial P^{\mathrm{ab}}=0\).  Thus
\(P^{\mathrm{ab}}\) is a \(\bbZ\)-circulation, equivalently an element of
\(H_1(G,\bbZ)\).  Dually, after applying the edge-coordinate pairing, it may be
viewed as an integral cohomology class; however, we keep the homological notation
used in the trace formulas.  Reducing coefficients modulo \(t\), we obtain
\[
   P^{\mathrm{ab}}_{[t]}\in H_1(G,\bbZ/t\bbZ).
\]
We say that \(P\) is \emph{nowhere-zero modulo \(t\)} if the coefficient of every
edge in \(P^{\mathrm{ab}}_{[t]}\) is nonzero in \(\bbZ/t\bbZ\).  We say that \(P\)
has \emph{integrally nowhere-zero abelianization} if the coefficient of every edge
in \(P^{\mathrm{ab}}\) is a nonzero integer.

Following the weighted setting of homological spectral graph
theory~\cite{LuoRoyHSGT}, let a \emph{positive directed weight} be a function
\[
   W:\mathbf E(G)\longrightarrow \bbR_{>0}.
\]
No symmetry between \(W(\mathbf e)\) and \(W(\mathbf e^{-1})\) is assumed.  For a
walk \(P=\mathbf e_{i_1}\cdots\mathbf e_{i_\ell}\), put
\[
   W(P)=\prod_{j=1}^{\ell}W(\mathbf e_{i_j}).
\]
Let \(A_W\) be the corresponding weighted vertex adjacency matrix, and define
\[
   \wtN_{W,\ell}(G)=\sum_{\substack{P\ \mathrm{closed}\\ |P|=\ell}}W(P),
   \qquad
   \wtN_{W,\ell}^{t\text{-}\NZ}(G)
      =\sum_{\substack{P\ \mathrm{closed},\ |P|=\ell\\
                        P^{\mathrm{ab}}_{[t]}\ \mathrm{nowhere\ zero}}}W(P).
\]
Thus \(\wtN_{W,\ell}(G)=\Tr(A_W^\ell)\).  Let \(B_W\) be the weighted
nonbacktracking edge adjacency matrix, and similarly put
\[
   N_{W,\ell}(G)=\Tr(B_W^\ell),
   \qquad
   N_{W,\ell}^{t\text{-}\NZ}(G)
      =\sum_{\substack{P\ \mathrm{circuit},\ |P|=\ell\\
                        P^{\mathrm{ab}}_{[t]}\ \mathrm{nowhere\ zero}}}W(P).
\]
For the integral condition, define \(\wtN_{W,\ell}^{\bbZ\text{-}\NZ}(G)\) and
\(N_{W,\ell}^{\bbZ\text{-}\NZ}(G)\) by the same weighted sums, with
\(P^{\mathrm{ab}}\) required to be nowhere zero over \(\bbZ\).  When \(W\equiv1\),
we omit \(W\) from the subscript; this convention applies to all six counting
functions above.  Finally, set
\[
   \rho_{A,W}=\rho(A_W),\qquad \rho_{B,W}=\rho(B_W),
\]
and let \(\nu_G\) be the period of the nonbacktracking edge shift, equivalently the
greatest common divisor of the lengths of all circuits in \(G\), in
the cases where this shift is nonempty.  The period depends only on \(G\), not on
the positive directed weight.  Suppose \(g\ge2\).  Following homological spectral
graph theory~\cite[Definition~3.5]{LuoRoyHSGT}, consider the assignment
\[
   P^{\mathrm{ab}}\longmapsto \frac{|P|}{\nu_G}\pmod 2
\]
on prime circuits.  If this assignment descends through their
integral homology relations to a homomorphism
\[
   \varepsilon_G:H_1(G,\bbZ)\longrightarrow\bbZ/2\bbZ,
\]
then
\[
   \eta(\alpha)=(-1)^{\varepsilon_G(\alpha)}
\]
is the \emph{period character} of \(G\).  It is unique whenever it exists.

\begin{definition}[Period modulating flow polynomial]
\label{def:period-modulating-flow-poly}
Let \(G\) be connected of genus \(g\ge2\), and suppose that its period character
\(\eta\) exists.  The \emph{modulating flow polynomial associated with the period
character} is
\begin{equation}\label{eq:period-modulating-flow-poly}
   \Deta(t)
      =
      \sum_{\substack{F\subseteq E\\
             \eta|_{H_1(G[F],\bbZ)}=1}}
         (-1)^{|E|-|F|}t^{\beta(F)}.
\end{equation}
\end{definition}

For later use, put
\[
   \Phi^A_{G,t}(\ell)
      =
      \begin{cases}
      \CG(t)+(-1)^\ell\DG(t),
         & G\ \text{is nonbipartite and }2\mid t,\\
      \CG(t),& \text{otherwise},
      \end{cases}
\]
and, when \(g\ge2\) and \(\nu_G\mid\ell\), put
\[
   \Phi^B_{G,t}(\ell)
      =
      \begin{cases}
      \CG(t)+(-1)^{\ell/\nu_G}\Deta(t),
         & 2\mid t\ \text{and the period character exists},\\
      \CG(t),& \text{otherwise}.
      \end{cases}
\]
Our main theorem is the following.

\begin{theorem}\label{thm:main}
Let \(G\) be a connected graph of genus \(g\), let
\(W:\mathbf E(G)\to\bbR_{>0}\) be a positive directed weight, and let \(t\ge 2\).
Then the following weighted asymptotic formulas hold as \(\ell\to\infty\).
\begin{enumerate}[label=(\alph*)]
\item For closed walks, there is a number \(\sigma_A<\rho_{A,W}\) such that
\begin{align*}
   \wtN_{W,\ell}^{t\text{-}\NZ}(G)
      &=
      \frac{\Phi^A_{G,t}(\ell)}{t^g}\,\wtN_{W,\ell}(G)
         +O(\sigma_A^\ell)\\
      &=
      \begin{cases}
      \dfrac{\Phi^A_{G,t}(\ell)}{t^g}\rho_{A,W}^\ell
         +O(\sigma_A^\ell),
         & G\ \text{is nonbipartite},\\[1ex]
      \dfrac{2\CG(t)}{t^g}\rho_{A,W}^\ell
         +O(\sigma_A^\ell),
         & G\ \text{bipartite},\ 2\mid \ell,\\[1ex]
      0, & G\ \text{bipartite},\ 2\nmid \ell.
      \end{cases}
\end{align*}
\item If \(g\ge2\), then there is a number \(\sigma_B<\rho_{B,W}\) such that,
for circuits,
\begin{align*}
   N_{W,\ell}^{t\text{-}\NZ}(G)
      &=\begin{cases}
      \dfrac{\Phi^B_{G,t}(\ell)}{t^g}N_{W,\ell}(G)+O(\sigma_B^\ell),
         & \nu_G\mid \ell,\\[1ex]
      0, & \nu_G\nmid \ell.
      \end{cases}\\
      &=\begin{cases}
      \dfrac{\nu_G\Phi^B_{G,t}(\ell)}{t^g}\rho_{B,W}^\ell
         +O(\sigma_B^\ell),
         & \nu_G\mid \ell,\\[1ex]
      0, & \nu_G\nmid \ell.
      \end{cases}
\end{align*}
\end{enumerate}
\end{theorem}

\begin{remark}[Why the bipartite case matters]
\label{rem:bipartite-edge-exception}
The canonical character is trivial on a bipartite graph, but this does not force the
trivial character to be the only edge spectral-radius maximizer.  The complete
classification in~\cite[Theorem~3.6(PC1)--(PC2), (SR4)--(SR5)]{LuoRoyHSGT} shows
that there is at most one further maximizer, namely the period character \(\eta\).
If \(G\) is nonbipartite, then \(\nu_G\) is odd, \(\eta=\theta\), and
\(\Deta=\DG\); hence \((-1)^{\ell/\nu_G}=(-1)^\ell\) along admissible lengths and
the vertex and edge corrections agree.  If \(G\) is bipartite, \(\eta\) may or may
not exist.  In particular, it does not exist for a regular bipartite graph of genus
at least two~\cite[Theorem~C.1(d)]{LuoRoyHSGT}; since existence is independent of
the weight, the absence of a correction then holds for every positive directed
weight.  Example~\ref{ex:bipartite-period-modulation} exhibits the other alternative.
\end{remark}

On every admissible subsequence on which the relevant total weighted count is
nonzero, Theorem~\ref{thm:main} implies convergence of the nowhere-zero ratio to
\(\Phi^A_{G,t}(\ell)/t^g\) or \(\Phi^B_{G,t}(\ell)/t^g\), respectively.  This
remains meaningful when the displayed coefficient is zero; only an asymptotic
equivalence with a nonzero leading constant would then be inappropriate.

Beyond the main asymptotic theorem, we investigate additional structural properties
of the character-indexed family.  For the full family,
Section~\ref{sec:matroidal-structure} realizes the retained edge sets as the balanced
sets of a complex gain graph and hence as the central sets of a semimatroid.  For the
canonical member \(D_G\) there is also an intrinsic matroidal description.  If
\(M=M(G)\) is the cycle matroid of \(G\), the classical cographic identity gives
\(\CG(t)=\chi_{M^*}(t)\)~\cite{Tutte1954}, while
Proposition~\ref{prop:matroidal-formulations} proves
\[
   \DG(t)=D_M(t),
\]
where \(D_M\) is the even-restriction polynomial obtained by retaining those
\(F\subseteq E\) for which every matroid circuit of \(M|F\) has even cardinality.
Thus \(\DG\) depends only on the cycle matroid.  The associated
Tutte specialization forces the coefficients of \(\DG(t)\), when nonzero, to
alternate strictly in sign with no internal zeros and gives the exact criterion
\[
   \DG(t)\not\equiv0
      \quad\Longleftrightarrow\quad
   G\ \text{is bridgeless}.
\]
Thus \(\DG\) has compatible Fourier, cycle-matroid, balanced-semimatroid, and
signed-Tutte descriptions; in the last it satisfies signed deletion-contraction.
The signed-Tutte framework also supplies the natural home for \(\Deta\), the
modulating flow polynomial associated with the switching class representing the
period character.

The rest of the paper explains why these two coefficients govern the leading terms.
Section~\ref{sec:fourier-coefficients} develops the general character coefficients,
specializes them to the trivial, canonical, and period characters, and derives their
parity consequences.  Section~\ref{sec:trace-selection} supplies the twisted trace
formulas and the complete spectral selection principle.
Section~\ref{sec:main-proof} proves Theorem~\ref{thm:main} by combining spectral
selection with edge inclusion-exclusion, derives the integral approximation theorem,
and applies it to weighted postman circuits.
Section~\ref{sec:matroidal-structure} establishes the gain-semimatroid structure of
the full family, the additional cycle-matroid structure of \(\DG\), and the
signed-Tutte framework for the order-two specializations.  Section~\ref{sec:examples}
gives exact and asymptotic examples of canonical and period modulation.

\section{Fourier coefficients and canonical--period modulation}
\label{sec:fourier-coefficients}

This section treats the nowhere-zero condition on the mod-\(t\) flow space by
inclusion-exclusion and character sums.  The resulting Fourier--Tutte expansion is
first packaged into a polynomial for every integral character.  Its trivial,
canonical, and period members are then identified with \(\CG\), \(\DG\), and
\(\Deta\), respectively.
We then derive the exact parity decompositions and fixed-modulus vanishing
obstructions.  The intrinsic matroidal, semimatroidal, and signed-Tutte structure is
developed separately in Section~\ref{sec:matroidal-structure}.

\subsection{The Fourier--Tutte expansion}

In this section \(G\) is an arbitrary finite graph.  Fix an orientation
\(\mathfrak o\) of \(G\).
For a positive integer \(t\), let
\[
   Z_t(G)=\ker\bigl(\partial:C_1(G,\bbZ/t\bbZ)\to C_0(G,\bbZ/t\bbZ)\bigr)
\]
be the mod-\(t\) flow space, and let
\[
   \mathcal F_t(G)
      =
      \{x\in Z_t(G):x_e\ne0\ \text{for all }e\in E\}
\]
be the set of nowhere-zero mod-\(t\) flows.  If \(F\subseteq E\), we write
\[
   Z_t(G[F])
      =
      \ker\bigl(\partial:C_1(G[F],\bbZ/t\bbZ)\to C_0(G,\bbZ/t\bbZ)\bigr),
\]
viewed as a subgroup of \(Z_t(G)\) by extending flows by zero outside \(F\).

For a character \(\xi:Z_t(G)\to\bbC^\ast\), define the \emph{Fourier coefficient}
\[
   \Psi_{G,\xi}(t)=\sum_{x\in\mathcal F_t(G)}\xi(x).
\]
The character-sum viewpoint below is in the same broad circle as the use of
finite-group, Ehrhart, group-arrangement, and Tutte-polynomial methods in the study
of flow
polynomials~\cite{Chen2007Orientations,ChenYang2009,BreuerSanyal2012,Chen2013DualComplementary}.

\begin{proposition}[Fourier--Tutte expansion]\label{prop:fourier-tutte}
Let \(G=(V,E)\) be a finite graph, let \(t\) be a positive integer, and let
\(\xi:Z_t(G)\to\bbC^\ast\) be a character.  Then
\[
   \Psi_{G,\xi}(t)
      =
      \sum_{\substack{F\subseteq E\\
         \xi|_{Z_t(G[F])}=1}}
      (-1)^{|E|-|F|}t^{\beta(F)}.
\]
\end{proposition}

\begin{proof}
Write \(\mathbf 1_A\) for the indicator of a condition \(A\).  For
\(x\in Z_t(G)\), the indicator that \(x\) is nowhere zero is
\[
   \mathbf 1_{\mathcal F_t(G)}(x)
      =
      \prod_{e\in E}(1-\mathbf 1_{\{x_e=0\}}).
\]
Expanding the product gives
\[
   \mathbf 1_{\mathcal F_t(G)}(x)
      =\sum_{S\subseteq E}(-1)^{|S|}
        \prod_{e\in S}\mathbf 1_{\{x_e=0\}}.
\]
After multiplying by \(\xi(x)\), summing over \(Z_t(G)\), and putting
\(F=E\setminus S\), we obtain
\[
\begin{aligned}
   \Psi_{G,\xi}(t)
      &=\sum_{x\in Z_t(G)}\xi(x)\mathbf 1_{\mathcal F_t(G)}(x)\\
      &=\sum_{F\subseteq E}(-1)^{|E|-|F|}
        \sum_{x\in Z_t(G)}\xi(x)
        \prod_{e\in E\setminus F}\mathbf 1_{\{x_e=0\}}\\
      &=\sum_{F\subseteq E}(-1)^{|E|-|F|}
        \sum_{x\in Z_t(G[F])}\xi(x).
\end{aligned}
\]
Indeed, the condition \(x_e=0\) for \(e\notin F\) says exactly that \(x\) is
supported on \(F\); restriction to \(F\), or equivalently extension by zero on
\(E\setminus F\), identifies these flows with \(Z_t(G[F])\).

We next evaluate the inner sum.  For a finite abelian group \(H\) and a character
\(\zeta:H\to\bbC^\ast\), character orthogonality gives
\[
   \sum_{x\in H}\zeta(x)
      =
      \begin{cases}
         |H|,&\zeta=1,\\
         0,&\zeta\ne1.
      \end{cases}
\]
For completeness, if \(\zeta\ne1\), choose \(y\in H\) with
\(\zeta(y)\ne1\).  Since translation by \(y\) permutes \(H\),
\[
   \sum_{x\in H}\zeta(x)
      =\sum_{x\in H}\zeta(x+y)
      =\zeta(y)\sum_{x\in H}\zeta(x).
\]
Because \(\zeta(y)\ne1\), this forces the character sum to vanish.  Applying this to
\(H=Z_t(G[F])\) and \(\zeta=\xi|_{Z_t(G[F])}\) shows that the inner sum vanishes
unless the restricted character is trivial.

It remains to record the order of the flow group.  Choose a spanning forest \(T\)
of the spanning subgraph \(G[F]\), with one tree in each of its \(c(F)\)
components.  Then
\[
   |T|=|V|-c(F),
   \qquad
   |F\setminus T|=|F|-|V|+c(F)=\beta(F).
\]
The fundamental cycles associated with the edges of \(F\setminus T\) form a basis
of the cycle module over \(\bbZ/t\bbZ\); this includes a loop as a one-edge
fundamental cycle.  Consequently,
\[
   Z_t(G[F])\cong(\bbZ/t\bbZ)^{\beta(F)},
   \qquad
   |Z_t(G[F])|=t^{\beta(F)}.
\]
Substitution in the inclusion--exclusion sum yields the asserted expansion.
\end{proof}

Proposition~\ref{prop:fourier-tutte} gives a character expansion separately for
each positive integer modulus \(t\), using a character of the finite flow group
\(Z_t(G)\).  We now show that a character of integral homology determines a
single polynomial whose value at every modulus to which the character descends
is the corresponding Fourier coefficient.  Let
\[
   \mathcal X_{\bbC}(G)
      =\operatorname{Hom}\bigl(H_1(G,\bbZ),\bbC^\ast\bigr)
\]
be the full complex character torus.  For a positive integer \(t\), let
\[
   r_t:H_1(G,\bbZ)\longrightarrow H_1(G,\bbZ/t\bbZ)
\]
be reduction modulo \(t\).  We say that \(\zeta\in\mathcal X_{\bbC}(G)\)
\emph{descends modulo \(t\)} if there is a character
\(\zeta_t:H_1(G,\bbZ/t\bbZ)\to\bbC^\ast\) such that
\(\zeta=\zeta_t\circ r_t\).  This is the specialization
\(\Lambda=tH_1(G,\bbZ)\) of the finite-quotient character formalism in
\cite[Lemma~2.13 and Theorem~3.21(b)]{LuoRoyHSGT} when \(G\) is connected; for
a general graph, the same construction is applied independently to its connected
components.

By \cite[Definition~2.11]{LuoRoyHSGT}, applied componentwise, every
\(\zeta\in\mathcal X_{\bbC}(G)\) is represented by a \(\bbC^\ast\)-gain graph
\(\Phi_\zeta\), unique up to switching, whose gain around an oriented cycle \(C\)
is \(\zeta([C])\), where \([C]\in H_1(G,\bbZ)\) is the class represented by
\(C\).  Since the oriented cycle classes of \(G[F]\) generate
\(H_1(G[F],\bbZ)\), an edge set \(F\) is balanced in \(\Phi_\zeta\) exactly when
\(\zeta|_{H_1(G[F],\bbZ)}=1\).  Define the \emph{character-indexed family}
\(\{D_{G,\zeta}\}_{\zeta\in\mathcal X_{\bbC}(G)}\) by
\begin{equation}\label{eq:general-character-polynomial}
   D_{G,\zeta}(q)
      =
      \sum_{\substack{F\subseteq E\\
         \zeta|_{H_1(G[F],\bbZ)}=1}}
      (-1)^{|E|-|F|}q^{\beta(F)}.
\end{equation}
The following lemma identifies each member as a one-variable specialization of
Zaslavsky's balanced dichromatic polynomial.
\begin{lemma}[Balanced-dichromatic specialization]
\label{lem:balanced-dichromatic-specialization}
Let \(\Phi_\zeta\) be any \(\bbC^\ast\)-gain-graph representative of the switching
class determined by \(\zeta\), and let \(Q^b_{\Phi_\zeta}(u,v)\) denote
Zaslavsky's balanced dichromatic polynomial~\cite[Section~3]{Zaslavsky1995BiasedIII}.
Then
\begin{equation}\label{eq:balanced-dichromatic-specialization}
   D_{G,\zeta}(q)=(-1)^{|E|+|V|}Q^b_{\Phi_\zeta}(-1,-q).
\end{equation}
In particular, the right-hand side is independent of the chosen representative,
and \(D_{G,\zeta}\in\bbZ[q]\).
\end{lemma}

\begin{proof}
Zaslavsky's subset expansion is
\[
   Q^b_{\Phi_\zeta}(u,v)
      =\sum_{F\in\mathcal B(\Phi_\zeta)}
        u^{b(F)}v^{|F|-|V|+b(F)},
\]
where \(b(F)\) is the number of balanced connected components of \(G[F]\).
For a balanced edge set, \(b(F)=c(F)\), so substitution of
\((u,v)=(-1,-q)\) gives \eqref{eq:general-character-polynomial} after
multiplication by \((-1)^{|E|+|V|}\), and shows that \(D_{G,\zeta}\in\bbZ[q]\).
Switching preserves balanced edge sets, so the right-hand side is independent of
the representative.
\end{proof}

Here \(q\) is formal.  Below we show that, whenever \(\zeta\) descends modulo
\(t\), \(D_{G,\zeta}(t)\) is the corresponding Fourier coefficient.  The
gain-semimatroid consequences are developed in
Proposition~\ref{prop:semimatroid-specialization}.

\begin{proposition}[General character coefficients]
\label{prop:integral-character-coefficient}
Let \(\zeta\in\mathcal X_{\bbC}(G)\) and let \(t\) be a positive integer.  Then
\(\zeta\) descends modulo \(t\) if and only if \(\zeta^t=1\).  In that case
the descended character \(\zeta_t\) is unique and
\begin{equation}\label{eq:general-character-fourier}
   D_{G,\zeta}(t)
      =\Psi_{G,\zeta_t}(t)
      =\sum_{x\in\mathcal F_t(G)}\zeta_t(x).
\end{equation}
\end{proposition}

\begin{proof}
Applying \cite[Lemma~2.5(c)]{LuoRoyHSGT} to each connected component of \(G\)
and taking direct sums, reduction induces the canonical isomorphism
\[
   H_1(G,\bbZ)/tH_1(G,\bbZ)
      \cong H_1(G,\bbZ/t\bbZ).
\]
Therefore \(\zeta\) factors through the mod-\(t\) group if and only if it is trivial
on \(tH_1(G,\bbZ)\), which is equivalent to
\(\zeta(t\alpha)=\zeta(\alpha)^t=1\) for every
\(\alpha\in H_1(G,\bbZ)\).  The quotient property also gives uniqueness of the
descent.

For every \(F\subseteq E\), the reduction map
\[
   H_1(G[F],\bbZ)\longrightarrow H_1(G[F],\bbZ/t\bbZ)
\]
is surjective.  Consequently, \(\zeta_t\) is trivial on \(Z_t(G[F])\) if and only
if \(\zeta\) is trivial on \(H_1(G[F],\bbZ)\).  Substitution in
Proposition~\ref{prop:fourier-tutte} gives \eqref{eq:general-character-fourier}.
\end{proof}

Pullback along reduction modulo \(t\) consequently identifies
\[
   \operatorname{Hom}\bigl(H_1(G,\bbZ/t\bbZ),\bbC^\ast\bigr)
      \cong \mathcal X_{\bbC}(G)[t],
   \qquad
   \mathcal X_{\bbC}(G)[t]=\{\zeta:\zeta^t=1\}.
\]
In particular, a character of finite order \(N\) descends precisely when
\(N\mid t\), whereas a character of infinite order never descends.

\newpage
\begin{corollary}[The flow polynomial as a Fourier coefficient]
\label{cor:flow-fourier-coefficient}
Let \(\mathbf 1\in\mathcal X_{\bbC}(G)\) and
\(\mathbf 1_t:Z_t(G)\to\bbC^\ast\) denote the integral and mod-\(t\) trivial
characters, respectively.  Then
\[
   \CG(t)
      =D_{G,\mathbf 1}(t)
      =\Psi_{G,\mathbf 1_t}(t)
      =\sum_{x\in\mathcal F_t(G)}\mathbf 1_t(x)
      =|\mathcal F_t(G)|.
\]
\end{corollary}

\begin{proof}
Every restriction of \(\mathbf 1\) is trivial, so
\eqref{eq:general-character-polynomial} is the defining expansion of \(\CG\).
The remaining identities follow from
Proposition~\ref{prop:integral-character-coefficient}.
\end{proof}

\subsection{The canonical coefficient}

For a connected graph \(K\), let
\(\theta_K:H_1(K,\bbZ)\to\{\pm1\}\) be the canonical character introduced in
homological spectral graph theory~\cite[Lemma~3.1 and
Definition~3.2]{LuoRoyHSGT}.  If \(G_1,\ldots,G_c\) are the connected components
of \(G\), use the canonical decomposition
\[
   H_1(G,\bbZ)=\bigoplus_{i=1}^c H_1(G_i,\bbZ)
\]
to define the \emph{canonical character} of an arbitrary graph \(G\) by
\begin{equation}\label{eq:disconnected-canonical-character}
   \theta_G(\alpha_1,\ldots,\alpha_c)
      =\prod_{i=1}^c\theta_{G_i}(\alpha_i).
\end{equation}
We write this character simply as \(\theta\).  Equivalently, with respect to the
orientation fixed above, if
\(\alpha=\sum_{\mathbf e\in\mathbf E_{\mathfrak o}(G)}
\alpha_{\mathbf e}\mathbf e\), then
\[
   \theta(\alpha)
      =(-1)^{\sum_{\mathbf e\in\mathbf E_{\mathfrak o}(G)}\alpha_{\mathbf e}}.
\]
This formula is independent of the orientation.  Applying
\cite[Theorem~3.6(SA2)--(SA3)]{LuoRoyHSGT} componentwise shows that \(\theta\)
takes values in \(\{\pm1\}\) and is trivial exactly when \(G\) is bipartite.

If \(t\) is even, then \(\theta^t=1\), since
\(\theta(t\alpha)=\theta(\alpha)^t=1\) for every
\(\alpha\in H_1(G,\bbZ)\).  Proposition~\ref{prop:integral-character-coefficient}
therefore gives the unique descent
\[
   \theta_t:H_1(G,\bbZ/t\bbZ)\longrightarrow\{\pm1\},
   \qquad
   \theta_t(\alpha\bmod t)=\theta(\alpha).
\]
With the orientation fixed above, the value of the descended character on
\(x\in Z_t(G)\) is
\[
   \theta_t(x)=(-1)^{\sum_{\mathbf e\in\mathbf E_{\mathfrak o}(G)}x_{\mathbf e}},
\]
where each coefficient \(x_{\mathbf e}\in\bbZ/t\bbZ\) is read modulo \(2\).  This expression is
independent of the chosen orientation.  If \(t\) is odd, then \(\theta^t=\theta\),
so Proposition~\ref{prop:integral-character-coefficient} shows that \(\theta\)
descends only when it is trivial.

\begin{proposition}[The canonical coefficient]\label{prop:canonical-coefficient}
Let \(t\) be a positive even integer, and let \(\theta_t\) be the mod-\(t\) descent
of the canonical character.  Then
\[
   \DG(t)
      =D_{G,\theta}(t)
      =\Psi_{G,\theta_t}(t)
      =\sum_{x\in\mathcal F_t(G)}\theta_t(x).
\]
If \(G_1,\ldots,G_c\) are the connected components of \(G\), then
\[
   \DG=\prod_{i=1}^c D_{G_i}.
\]
Moreover, let \(c_{\mathrm{nb}}(G)\) be the number of nonbipartite connected
components of \(G\).  If \(G\) is bipartite, then \(\DG=\CG\).  If
\(\DG\not\equiv0\), then
\[
   \deg\DG\le \beta(G)-c_{\mathrm{nb}}(G),
\]
with equality \(\deg\DG=\beta(G)\) when \(G\) is bipartite.
\end{proposition}

\begin{proof}
By the definition of \(\Psi_{G,\theta_t}(t)\), it remains to identify the subgraphs
\(G[F]\) on whose flow space \(\theta_t\) is trivial.  For an oriented cycle \(C\),
viewed as a mod-\(t\) flow, the descended canonical character evaluates as
\[
   \theta_t(C)=(-1)^{|C|}.
\]
Hence \(\theta_t\) is nontrivial on \(Z_t(G[F])\) whenever \(G[F]\) contains an odd
cycle.  Conversely, if \(G[F]\) is bipartite, every cycle has even length, and
reducing flows modulo \(2\) shows that every element of \(Z_t(G[F])\) has even total
coordinate sum.  Thus \(\theta\) is trivial on \(H_1(G[F],\bbZ)\), equivalently
\(\theta_t\) is trivial on \(Z_t(G[F])\), exactly when \(G[F]\) is bipartite.
Hence \(D_{G,\theta}=D_G\), and
Proposition~\ref{prop:integral-character-coefficient} gives the displayed Fourier
identity.

To prove the component factorization, let \(F\subseteq E\).  Then \(F\)
decomposes uniquely as \(F=F_1\sqcup\cdots\sqcup F_c\), where
\(F_i\subseteq E(G_i)\), and
\[
   \beta(F)=\sum_{i=1}^c\beta(F_i).
\]
Moreover, \(G[F]\) is bipartite if and only if every \(G_i[F_i]\) is bipartite.
The defining subset expansion of \(\DG\) therefore factors as asserted.

It remains to establish the degree assertions for a connected graph.  If \(G\)
is bipartite, then every spanning
subgraph \(G[F]\) is bipartite, so the defining expansions of \(\DG\) and \(\CG\)
coincide.  If \(G\) has a bridge, both polynomials vanish; otherwise the term
\(F=E\) is the unique term of cycle rank \(\beta(G)\), and hence
\(\deg\DG=\beta(G)\).

If \(G\) is nonbipartite, consider a bipartite spanning subgraph \(G[F]\).  Were
\(\beta(F)=\beta(G)\), the identity
\[
   \beta(G)-\beta(F)=|E\setminus F|-c(F)+1
\]
would imply that the omitted edges only connect the components of \(G[F]\) and
create no additional cycle.  Thus
\(H_1(G[F],\bbZ)=H_1(G,\bbZ)\).  This is impossible: the canonical character is
trivial on the former group because \(G[F]\) is bipartite, but nontrivial on the
latter because \(G\) is nonbipartite.  Consequently every term in \(\DG\) has
cycle rank at most \(\beta(G)-1\).

Finally, apply these connected conclusions to the component factors.  If none
of the factors vanishes, their
degrees add.  Each nonbipartite component lowers the resulting upper bound by
at least one, giving
\(\deg\DG\le \beta(G)-c_{\mathrm{nb}}(G)\).
\end{proof}

\subsection{The period coefficient}

For a connected graph \(K\) of genus at least two, let \(\eta_K\) denote the
period character of~\cite[Definition~3.5]{LuoRoyHSGT}, when it exists.  Suppose
that the connected components \(G_1,\ldots,G_c\) of \(G\) all have genus at least
two and admit period characters \(\eta_{G_1},\ldots,\eta_{G_c}\).  In analogy
with \eqref{eq:disconnected-canonical-character}, define the \emph{componentwise
period character} by
\begin{equation}\label{eq:componentwise-period-character}
   \eta_G(\alpha_1,\ldots,\alpha_c)
      =\prod_{i=1}^c\eta_{G_i}(\alpha_i),
   \qquad
   (\alpha_1,\ldots,\alpha_c)\in H_1(G,\bbZ).
\end{equation}
We write this character simply as \(\eta\).  The qualifier ``componentwise'' is
essential: the nonbacktracking matrix of a disconnected graph is reducible, and
the periods \(\nu_{G_i}\) of its connected blocks need not agree.  Thus
\eqref{eq:componentwise-period-character} defines a homology character, not a
single Perron--Frobenius period for the whole graph.

By \cite[Theorem~3.6(PC1)]{LuoRoyHSGT}, every \(\eta_{G_i}\) is a nontrivial
character of order two.  Hence \(\eta\) has order two.  If \(t\) is even, then
\(\eta^t=1\), and Proposition~\ref{prop:integral-character-coefficient} gives the
unique descent
\[
   \eta_t:H_1(G,\bbZ/t\bbZ)\longrightarrow\{\pm1\},
   \qquad
   \eta_t(\alpha\bmod t)=\eta(\alpha),
\]
For odd \(t\), \(\eta^t=\eta\ne1\), so the same proposition shows that \(\eta\)
does not descend modulo \(t\).

Thus the polynomial in \eqref{eq:period-modulating-flow-poly} is exactly the
period member \(D_{G,\eta}\) of the character-indexed family and the corresponding
modulating flow polynomial.

\begin{proposition}[The period coefficient]\label{prop:period-coefficient}
Let \(t\) be a positive even integer, and suppose that the componentwise period
character \(\eta\) of \(G\) exists.  Then
\[
   \Deta(t)
      =\Psi_{G,\eta_t}(t)
      =\sum_{x\in\mathcal F_t(G)}\eta_t(x).
\]
Moreover, if \(G_1,\ldots,G_c\) are the connected components of \(G\), then
\[
   D_{G,\eta}=\prod_{i=1}^cD_{G_i,\eta_{G_i}},
\]
and, if \(D_{G,\eta}\not\equiv0\), then
\[
   \deg D_{G,\eta}\le\beta(G)-c.
\]
If every
connected component of \(G\) is nonbipartite, then \(\eta=\theta\), hence
\(\eta_t=\theta_t\), and
\[
   \Deta(t)=\DG(t).
\]
\end{proposition}

\begin{proof}
The Fourier identity follows directly from
Proposition~\ref{prop:integral-character-coefficient}, since \(\eta\) has order two.

Write \(F=F_1\sqcup\cdots\sqcup F_c\), where \(F_i\subseteq E(G_i)\).  Under the
corresponding direct-sum decomposition of \(H_1(G[F],\bbZ)\), the restriction of
\(\eta\) is trivial if and only if every restriction
\(\eta_{G_i}|_{H_1(G_i[F_i],\bbZ)}\) is trivial.  Since cycle ranks and signs also
split over the components, the defining subset expansion gives the displayed
factorization.  Each period character \(\eta_{G_i}\) is nontrivial by
\cite[Theorem~3.6(PC1)]{LuoRoyHSGT}.  The connected cycle-rank argument in the
proof of Proposition~\ref{prop:canonical-coefficient}, with \(\eta_{G_i}\) in
place of the canonical character, therefore gives
\(\deg D_{G_i,\eta_{G_i}}\le\beta(G_i)-1\) whenever that factor is nonzero.
Adding these bounds proves \(\deg D_{G,\eta}\le\beta(G)-c\).

Finally, if every \(G_i\) is nonbipartite, then
\(\eta_{G_i}=\theta_{G_i}\) by \cite[Theorem~3.6(PC1)]{LuoRoyHSGT}.  Equations
\eqref{eq:disconnected-canonical-character} and
\eqref{eq:componentwise-period-character} therefore give \(\eta=\theta\), and
Proposition~\ref{prop:canonical-coefficient} gives \(D_{G,\eta}=D_G\).
\end{proof}

\subsection{Parity bias and flow obstructions}

Every order-two character gives an exact two-class bias inside the nowhere-zero
flow set.  The canonical and period characters are the cases selected later by the
closed-walk and nonbacktracking dynamics, respectively.

\begin{proposition}[Parity bias]\label{prop:parity-bias}
Let \(G\) be a finite graph, let
\(\zeta\in\mathcal X_{\bbC}(G)\) satisfy \(\zeta^2=1\), and let \(t\) be a positive
even integer.  Then \(\zeta\) has a unique mod-\(t\) descent \(\zeta_t\), which is
\(\{\pm1\}\)-valued.  Put
\[
   \mathcal F_{t,\zeta}^+(G)=\{x\in\mathcal F_t(G):\zeta_t(x)=1\},
   \qquad
   \mathcal F_{t,\zeta}^-(G)=\{x\in\mathcal F_t(G):\zeta_t(x)=-1\}.
\]
Then
\[
   |\mathcal F_{t,\zeta}^+(G)|=\frac{\CG(t)+D_{G,\zeta}(t)}{2},
   \qquad
   |\mathcal F_{t,\zeta}^-(G)|=\frac{\CG(t)-D_{G,\zeta}(t)}{2}.
\]
In particular,
\[
   |D_{G,\zeta}(t)|\le \CG(t),
   \qquad
   \CG(t)\equiv D_{G,\zeta}(t)\pmod 2.
\]
\end{proposition}

\begin{proof}
Since \(t\) is even and \(\zeta^2=1\), we have \(\zeta^t=1\).  Hence
Proposition~\ref{prop:integral-character-coefficient} gives both the unique descent
\(\zeta_t\) and the Fourier identity
\[
   D_{G,\zeta}(t)=\sum_{x\in\mathcal F_t(G)}\zeta_t(x).
\]
Since reduction modulo \(t\) is surjective, \(\zeta_t^2=1\), so the two sets in the
statement partition \(\mathcal F_t(G)\).
The ordinary flow polynomial counts nowhere-zero mod-\(t\) flows, so
\[
   \CG(t)=|\mathcal F_t(G)|
      =|\mathcal F_{t,\zeta}^+(G)|+|\mathcal F_{t,\zeta}^-(G)|.
\]
Moreover, the Fourier identity gives
\[
   D_{G,\zeta}(t)
      =|\mathcal F_{t,\zeta}^+(G)|-|\mathcal F_{t,\zeta}^-(G)|.
\]
Solving these two linear equations gives the two displayed formulas.  The inequality
and congruence are immediate.
\end{proof}

\begin{corollary}[Flow obstruction]\label{cor:flow-obstruction}
Let \(G\) be a finite graph and let \(t\) be a positive even integer.  If
\(\CG(t)=0\), then \(\DG(t)=0\) and, whenever the componentwise period character
of \(G\) exists, \(\Deta(t)=0\).
\end{corollary}

\begin{proof}
Both assertions follow from Proposition~\ref{prop:parity-bias}.
\end{proof}

\begin{corollary}\label{cor:noneulerian}
If \(G\) is not Eulerian, then \(\DG(2)=0\).
\end{corollary}

\begin{proof}
The classical characterization of nowhere-zero \(2\)-flows states that
\(\CG(2)>0\) if and only if every vertex of \(G\) has even degree; see, for example,
\cite{Zhang1997}.  Hence \(\CG(2)=0\), and
Corollary~\ref{cor:flow-obstruction} gives \(\DG(2)=0\).
\end{proof}

\paragraph{Relation with classical flow conjectures}
Tutte's classical flow conjectures concern nonvanishing of the ordinary flow
polynomial~\cite{Tutte1954,Zhang1997,Shahmohamad2013,Dong2020FlowZeros}.  The
\(5\)-flow conjecture asserts that \(\CG(5)>0\) for every bridgeless graph.  On a
nonbipartite graph the canonical character is nontrivial \(2\)-torsion and does not
descend to the mod-\(5\) flow space; on a bipartite graph it is already trivial.
Thus no additional canonical coefficient occurs in the corresponding case of
Theorem~\ref{thm:main}.  The Petersen-minor form of the \(4\)-flow
conjecture asserts that \(\CG(4)>0\) for every bridgeless graph with no Petersen
minor.  Here \(\DG(4)\) refines existence by measuring the imbalance
\[
   |\mathcal F_{4,\theta_4}^+(G)|
      =\frac{\CG(4)+\DG(4)}{2},
   \qquad
   |\mathcal F_{4,\theta_4}^-(G)|
      =\frac{\CG(4)-\DG(4)}{2}.
\]
If \(G\) has no nowhere-zero \(4\)-flow, Corollary~\ref{cor:flow-obstruction}
forces \(\DG(4)=0\); when such flows exist, \(\DG(4)\) measures their canonical
parity bias rather than providing a different existence criterion.

\section{Trace formulas and spectral selection}
\label{sec:trace-selection}

This section provides the spectral bridge between the Fourier coefficients of the
finite flow space and weighted walk counting.  We first twist the weighted vertex
and nonbacktracking edge adjacency matrices by \(\mu_t\)-valued gains and use Fourier
inversion to obtain exact trace formulas for prescribed mod-\(t\) abelianization.  We
then determine which character classes can contribute at spectral radius: the
trivial and canonical classes for closed walks, and the trivial class together with
the period class, when the latter exists, for circuits.
Section~\ref{sec:main-proof} combines this spectral selection with edge
inclusion-exclusion to prove the main asymptotic theorem.  Except in
Proposition~\ref{prop:spectral-input}, no connectivity assumption is imposed
anywhere in this section.

\subsection{Weighted twisted trace formulas}

To connect the nowhere-zero counting problem with the spectral estimates below, we
need exact trace identities that retain the positive directed weight.  Indeed, edge
inclusion--exclusion reduces the nowhere-zero condition to prescribed mod-\(t\)
coordinates on subsets of edges, and Fourier inversion must express each resulting
weighted count as a sum of traces of the corresponding twisted weighted transfer
operators.  Thus the weighted supported-coordinate formulation is an essential
input to the main asymptotic theorem.

When \(W\equiv1\), the supported-coordinate counterparts of
Propositions~\ref{prop:trace-formula}
and~\ref{prop:edge-trace-formula} are proved in
\cite[Theorem~3.3]{LuoEulerian}, with the homology specialization given in
\cite[Theorem~3.4]{LuoEulerian}.  The nonbacktracking trace formula in its original
finite-quotient homology form appears in
\cite[Theorem~3.21(b)]{LuoRoyHSGT}.  We establish here the precise weighted versions
needed for arbitrary positive directed weights.

Fix an orientation \(\mathfrak o\) of \(G\).  We identify \(C_1(G,\bbZ/t\bbZ)\) with
the free \(\bbZ/t\bbZ\)-module on \(\mathbf E_{\mathfrak o}(G)\).  Put
\[
   \mu_t=\{z\in\bbC^*:z^t=1\}.
\]

Following \cite[Definitions~2.11 and 2.18]{LuoRoyHSGT}, a \emph{complex gain} on
\(G\) is a function
\[
   \gamma:\mathbf E(G)\longrightarrow\bbC^*,
   \qquad
   \gamma(\mathbf e^{-1})=\gamma(\mathbf e)^{-1}.
\]
For a walk \(P=\mathbf e_{i_1}\cdots\mathbf e_{i_\ell}\), put
\[
   \gamma(P)=\prod_{j=1}^{\ell}\gamma(\mathbf e_{i_j}).
\]

Fix a positive directed weight \(W:\mathbf E(G)\to\bbR_{>0}\).  The trivial
specialization \(W\equiv1\) will always recover the unweighted objects.

\begin{definition}\label{def:twisted-vertex}
The \emph{twisted weighted vertex adjacency matrix} \(A_{\gamma,W}\) is the
\(|V|\times |V|\) matrix with entries
\[
   (A_{\gamma,W})_{uv}
      =
      \sum_{\substack{\mathbf e\in \mathbf E(G)\\
             \mathbf e(0)=u,\ \mathbf e(1)=v}}
        W(\mathbf e)\gamma(\mathbf e).
\]
When \(\gamma\equiv1\), we write \(A_W\), and when \(W\equiv1\), we write
\(A_\gamma=A_{\gamma,1}\).
\end{definition}

\begin{definition}\label{def:twisted-edge}
The \emph{twisted weighted nonbacktracking edge adjacency matrix}
\(B_{\gamma,W}\) is indexed by \(\mathbf E(G)\) and has entries
\[
   (B_{\gamma,W})_{\mathbf e\mathbf f}
      =
      \begin{cases}
      W(\mathbf f)\gamma(\mathbf f),
         & \mathbf e(1)=\mathbf f(0)\ \mathrm{and}\ \mathbf f\ne\mathbf e^{-1},\\
      0, & \mathrm{otherwise}.
      \end{cases}
\]
When \(\gamma\equiv1\), we write \(B_W\), and when \(W\equiv1\), we write
\(B_\gamma=B_{\gamma,1}\); in particular, \(B_0\) is the ordinary
nonbacktracking edge adjacency matrix.
\end{definition}

If \(S\subseteq E\), let \(M_t(S)\) be the free \(\bbZ/t\bbZ\)-module on the edges in
\(S\), viewed as the subgroup of \(C_1(G,\bbZ/t\bbZ)\) supported on \(S\).
Let \(\mathcal G_t(S)\) be the group of \(\mu_t\)-valued gains that are trivial on
every edge outside \(S\).  Restriction to the positively oriented edges over \(S\)
identifies \(\mathcal G_t(S)\) with \(\mu_t^S\), so
\[
   |\mathcal G_t(S)|=t^{|S|}.
\]
If \(\alpha=\sum_{e\in S}\alpha_e\mathbf e\in M_t(S)\) and
\(\gamma\in\mathcal G_t(S)\), define
\[
   \gamma(\alpha)=\prod_{e\in S}\gamma(\mathbf e)^{\alpha_e}.
\]
This is well defined because every value of \(\gamma\) is a \(t\)-th root of unity.
For \(\alpha\in M_t(S)\), let
\[
   \wtN_{W,t,S}(\alpha,\ell)
\]
be the sum of \(W(P)\) over the closed walks \(P\) of length \(\ell\) such that
\[
   P^{\mathrm{ab}}_{[t]}\big|_S=\alpha.
\]

\begin{proposition}[Weighted Fourier trace formula]\label{prop:trace-formula}
For every \(S\subseteq E\), every \(\alpha\in M_t(S)\), and every \(\ell\ge 1\),
\[
   \wtN_{W,t,S}(\alpha,\ell)
      =
      \frac{1}{t^{|S|}}
      \sum_{\gamma\in\mathcal G_t(S)}
         \gamma(\alpha)^{-1}
         \Tr(A_{\gamma,W}^\ell).
\]
\end{proposition}

\begin{proof}
Write \(A=A_{\gamma,W}\).  Expanding the diagonal entries of \(A^\ell\) first
over vertices gives
\[
   \Tr(A^\ell)
      =
      \sum_{v_0,\ldots,v_{\ell-1}\in V}
         A_{v_0v_1}A_{v_1v_2}\cdots A_{v_{\ell-1}v_0}.
\]
For each factor \(A_{v_{j-1}v_j}\), its definition is a sum over all oriented
edges from \(v_{j-1}\) to \(v_j\).  Distributing these edge sums therefore
identifies every term with a unique rooted closed walk
\[
   P=\mathbf e_{i_1}\cdots\mathbf e_{i_\ell},
   \qquad
   \mathbf e_{i_j}(1)=\mathbf e_{i_{j+1}}(0),
   \qquad
   \mathbf e_{i_\ell}(1)=\mathbf e_{i_1}(0).
\]
This expansion remains valid for loops and parallel edges because the choice of an
oriented edge is retained when the corresponding matrix entry is distributed.
The contribution of \(P\) is
\[
   \prod_{j=1}^{\ell}
      W(\mathbf e_{i_j})\gamma(\mathbf e_{i_j})
      =
      W(P)\gamma(P).
\]
Indeed, the weight factors multiply to \(W(P)\), while multiplicativity of
the gain turns the remaining factors into \(\gamma(P)\).  Since
\(\gamma\in\mathcal G_t(S)\), this value depends only on the mod-\(t\)
abelianization and equals \(\gamma(P^{\mathrm{ab}}_{[t]}|_S)\).  Consequently,
\begin{equation}\label{eq:weighted-vertex-trace-expansion}
   \Tr(A_{\gamma,W}^\ell)
      =
      \sum_{\substack{P\ \mathrm{closed}\\ |P|=\ell}}
         W(P)\gamma\!\left(P^{\mathrm{ab}}_{[t]}\big|_S\right).
\end{equation}
Substituting \eqref{eq:weighted-vertex-trace-expansion} into the Fourier sum and
interchanging the two finite sums yields
\[
 \frac{1}{t^{|S|}}\sum_{\gamma\in\mathcal G_t(S)}
      \gamma(\alpha)^{-1}\Tr(A_{\gamma,W}^\ell)
   =
   \sum_{\substack{P\ \mathrm{closed}\\ |P|=\ell}}W(P)
   \frac{1}{t^{|S|}}\sum_{\gamma\in\mathcal G_t(S)}
      \gamma\!\left(P^{\mathrm{ab}}_{[t]}|_S-\alpha\right).
\]
To evaluate the inner sum, put
\[
   \beta=P^{\mathrm{ab}}_{[t]}\big|_S-\alpha\in M_t(S).
\]
Using the oriented edges over \(S\) as a basis, write
\(\beta=(\beta_e)_{e\in S}\).  Then the
orthogonality sum factors coordinatewise:
\[
   \frac{1}{t^{|S|}}
      \sum_{\gamma\in\mathcal G_t(S)}\gamma(\beta)
      =
      \prod_{e\in S}
      \left(
         \frac1t\sum_{z\in\mu_t}z^{\beta_e}
      \right)
      =
      \prod_{e\in S}\mathbf 1_{\{\beta_e=0\}}
      =
      \mathbf 1_{\{\beta=0\}}.
\]
Thus a closed walk contributes \(W(P)\) precisely when
\(P^{\mathrm{ab}}_{[t]}|_S=\alpha\), and contributes \(0\) otherwise.  The outer
sum is therefore exactly \(\wtN_{W,t,S}(\alpha,\ell)\), proving the proposition.
\end{proof}

For \(\alpha\in M_t(S)\), let \(N_{W,t,S}(\alpha,\ell)\) be the sum of \(W(P)\)
over length-\(\ell\) circuits satisfying
\(P^{\mathrm{ab}}_{[t]}|_S=\alpha\).

\begin{proposition}[Weighted edge trace formula]\label{prop:edge-trace-formula}
For every \(S\subseteq E\), every \(\alpha\in M_t(S)\), and every \(\ell\ge1\),
\[
   N_{W,t,S}(\alpha,\ell)
      =
      \frac{1}{t^{|S|}}
      \sum_{\gamma\in\mathcal G_t(S)}
         \gamma(\alpha)^{-1}
         \Tr(B_{\gamma,W}^\ell).
\]
\end{proposition}

\begin{proof}
Expanding \(\Tr(B_{\gamma,W}^\ell)\) gives
\[
   \Tr(B_{\gamma,W}^\ell)
      =
      \sum_{\substack{P\ \mathrm{circuit}\\ |P|=\ell}}
         W(P)\gamma(P).
\]
Applying the same character-orthogonality calculation as in
Proposition~\ref{prop:trace-formula} selects exactly the circuits satisfying
\(P^{\mathrm{ab}}_{[t]}|_S=\alpha\).
\end{proof}

\subsection{Trivial, canonical, and period character classes}

A \(\mu_t\)-valued gain \(\gamma\) determines a character of
\(H_1(G,\bbZ/t\bbZ)\) by evaluating the gain on closed walks.  Two such gains
determine the same character exactly when they are switching equivalent.
Switching-equivalent gains give similar twisted weighted vertex and edge adjacency
matrices~\cite[Definition~2.11 and Lemma~2.19]{LuoRoyHSGT}.  Thus the spectra of
\(A_{\gamma,W}\) and \(B_{\gamma,W}\) depend only on the character represented by
\(\gamma\).

The trivial character class consists of the gains switching equivalent to the
trivial gain.  When \(t\) is even, the canonical character \(\theta\) descends to
\(\theta_t\), and the gains representing \(\theta_t\) form the canonical character
class; see
\cite[Lemma~3.1 and Definition~3.2]{LuoRoyHSGT}.  For every closed walk \(P\),
\[
   \theta_t(P^{\mathrm{ab}}_{[t]})=(-1)^\ell,
\]
where \(\ell\) is the length of \(P\).  This is the mod-\(t\) avatar of the canonical
character in homological spectral graph theory.  When \(t\) is odd, a nontrivial
\(\theta\) does not descend, so there is no nontrivial canonical class.

The following counting lemma applies simultaneously to the canonical and period
classes.  It is the point at which the two modulating polynomials enter the trace
formulas.

\begin{lemma}\label{lem:supported-character-counts}
Let \(G\) be an arbitrary finite graph, let \(F\subseteq E\), and put
\(S=E\setminus F\).
\begin{enumerate}[label=(\alph*)]
\item The number of gains \(\gamma\in\mathcal G_t(S)\) representing the trivial
      character on \(H_1(G,\bbZ/t\bbZ)\) is \(t^{c(F)-c(G)}\).
\item Suppose that \(t\) is even and let
      \(\zeta:H_1(G,\bbZ)\to\{\pm1\}\) be a nontrivial order-two character, with
      mod-\(t\) reduction \(\zeta_t\).  The number of gains
      \(\gamma\in\mathcal G_t(S)\) representing \(\zeta_t\) is
      \[
         \begin{cases}
         t^{c(F)-c(G)},&
            \zeta|_{H_1(G[F],\bbZ)}=1,\\
         0,&\text{otherwise}.
         \end{cases}
      \]
\end{enumerate}
\end{lemma}

\begin{proof}
A gain represents the trivial character if and only if it is switching equivalent
to the trivial gain.  Thus it has the form
\[
   \gamma_s(\mathbf e)=s(\mathbf e(0))^{-1}s(\mathbf e(1)),
   \qquad s:V\longrightarrow\mu_t.
\]
The condition \(\gamma_s\in\mathcal G_t(S)\) says that
\(\gamma_s(\mathbf e)=1\) for every edge in \(F\), or equivalently that \(s\) is
constant on each connected component of \(G[F]\).  There are \(t^{c(F)}\) such
switching functions.  Two switching functions give the same gain exactly when
their quotient is constant independently on each connected component of \(G\),
which gives a redundancy of \(t^{c(G)}\).  Hence there are
\(t^{c(F)-c(G)}\) gains in \(\mathcal G_t(S)\) representing the trivial character.

Now assume that \(t\) is even, and choose a \(\mu_t\)-valued gain
\(\gamma_\zeta\) representing \(\zeta_t\).  Every gain representing \(\zeta_t\)
is a switching of \(\gamma_\zeta\), of the form
\[
   \gamma_\zeta^s(\mathbf e)
      =s(\mathbf e(0))^{-1}\gamma_\zeta(\mathbf e)s(\mathbf e(1)).
\]
Requiring \(\gamma_\zeta^s\in\mathcal G_t(S)\) means making the switched gain
trivial on \(G[F]\).  This is possible exactly when the gain of
\(\gamma_\zeta\) around every cycle of \(G[F]\) is trivial, or equivalently when
\(\zeta|_{H_1(G[F],\bbZ)}=1\).  When it is possible, the solutions form a torsor
for the switching functions constant on each component of \(G[F]\).  Quotienting
their number \(t^{c(F)}\) by the \(t^{c(G)}\) switching functions constant on the
components of \(G\) gives \(t^{c(F)-c(G)}\) distinct gains.
This proves (b).  Taking \(\zeta=\theta\) recovers
the bipartite criterion in the former version of the lemma; taking \(\zeta=\eta\)
gives the support count for \(\Deta\).
\end{proof}

We shall use the following spectral input.  The full-character weighted spectral
antisymmetry theorem applies to arbitrary positive directed weights and pairs every
character with its canonical translate.  Its unitary spectral-radius extremality
statements identify the possible leading character classes for both vertex and
nonbacktracking edge matrices.  The edge classification is governed by the period
character, not by an
unspecified collection of order-two exceptions.  Write
\(\ee(x)=\exp(2\pi i x)\).  We use the following consequences.

\begin{proposition}[Spectral input from homological spectral graph
theory~\cite{LuoRoyHSGT}]\label{prop:spectral-input}
Let \(G\) be connected, let \(W\) be a positive directed weight, and put
\[
   \rho_{A,W}=\rho(A_W),\qquad \rho_{B,W}=\rho(B_W).
\]
For \(\mu_t\)-valued gains, the following statements hold.
\begin{enumerate}[label=(\alph*)]
\item If \(\gamma\) represents the trivial character, then
      \[
         \Tr(A_{\gamma,W}^\ell)=\wtN_{W,\ell}(G),
         \qquad
         \Tr(B_{\gamma,W}^\ell)=N_{W,\ell}(G).
      \]
\item If \(t\) is even and \(\gamma\) represents \(\theta_t\), then
      \[
         \Tr(A_{\gamma,W}^\ell)=(-1)^\ell\wtN_{W,\ell}(G),
         \qquad
         \Tr(B_{\gamma,W}^\ell)=(-1)^\ell N_{W,\ell}(G).
      \]
\item Suppose \(g\ge2\), the period character exists, and \(t\) is even.  If
      \(\gamma\) represents \(\eta_t\), then
      \[
         \Tr(B_{\gamma,W}^\ell)
            =\ee\!\left(\frac{\ell}{2\nu_G}\right)N_{W,\ell}(G).
      \]
      In particular, when \(\nu_G\mid\ell\), this factor is
      \((-1)^{\ell/\nu_G}\).
\item If \(\gamma\) represents neither the trivial character nor, when it is
      present, the canonical character, then
      \[
         \rho(A_{\gamma,W})<\rho_{A,W}.
      \]
\item Suppose \(g\ge2\).  If \(\gamma\) represents neither the trivial character
      nor, when \(t\) is even and \(\eta\) exists, the period character, then
      \[
         \rho(B_{\gamma,W})<\rho_{B,W}.
      \]
\end{enumerate}
\end{proposition}

\begin{proof}
Part (a) is switching equivalence, and part (b) is canonical spectral
antisymmetry~\cite[Theorem~3.6(SA1)]{LuoRoyHSGT}.  Period spectral rotation gives
\[
   \Spec B_{\gamma_\eta\gamma,W}
      =\ee\!\left(\frac{1}{2\nu_G}\right)\Spec B_{\gamma,W}
\]
whenever \(\gamma_\eta\) represents the period character and \(\gamma\) is a
complex gain, with the product of gains taken edgewise.  Taking
\(\gamma\equiv1\), raising eigenvalues to the \(\ell\)-th power, and summing
proves (c).  The vertex strictness in (d) is
\cite[Theorem~3.6(SR2)]{LuoRoyHSGT}.  Finally,
\cite[Theorem~3.6(SR4)--(SR5)]{LuoRoyHSGT} states that, for the edge matrices, the
unitary equality locus consists of the trivial and period characters when the
latter exists, and only the trivial character otherwise.  A nontrivial order-two
character belongs to the finite mod-\(t\) dual exactly when \(t\) is even, which
gives (e).
\end{proof}

\section{Proof of the main theorem}
\label{sec:main-proof}

This section has three parts.  First we prove the fixed-modulus asymptotic theorem by
combining edge inclusion-exclusion with the spectral selection principle from
Section~\ref{sec:trace-selection}.  We then let the modulus grow and use the degree
estimates from Section~\ref{sec:fourier-coefficients} to pass from mod-\(t\)
nowhere-zero abelianization to integral
nowhere-zero abelianization.  Finally, we apply the integral approximation to
weighted postman circuits.

\subsection{Fixed-modulus asymptotics}

\begin{proof}[Proof of Theorem~\ref{thm:main}]
We first prove the closed-walk statement.  For a closed walk \(P\), write
\[
   z_e(P)=\mathbf 1_{\{P^{\mathrm{ab}}_{[t]}(e)=0\}}
   \qquad (e\in E).
\]
Then \(P\) is nowhere-zero modulo \(t\) precisely when \(z_e(P)=0\) for every edge.
Thus
\[
   \mathbf 1_{\mathrm{nowhere\ zero}}(P)
      =
      \prod_{e\in E}(1-z_e(P))
      =
      \sum_{S\subseteq E}(-1)^{|S|}
      \mathbf 1_{\{P^{\mathrm{ab}}_{[t]}|_S=0\}}.
\]
Multiplying by \(W(P)\), summing over all length-\(\ell\) closed walks, and applying
Proposition~\ref{prop:trace-formula} with \(\alpha=0\), we obtain
\begin{align}
   \wtN_{W,\ell}^{t\text{-}\NZ}(G)
      &=
      \sum_{F\subseteq E}(-1)^{|E|-|F|}
      \frac{1}{t^{|E|-|F|}}
      \sum_{\gamma\in\mathcal G_t(E\setminus F)}\Tr(A_{\gamma,W}^\ell).
      \label{eq:ie-trace}
\end{align}
Here we have written \(S=E\setminus F\).

We now isolate the terms in \eqref{eq:ie-trace} which can have spectral radius
\(\rho_{A,W}\).  First consider gains in \(\mathcal G_t(E\setminus F)\) which
represent the trivial character.  By Proposition~\ref{prop:spectral-input}(a),
every such gain has
\[
   \Tr(A_{\gamma,W}^\ell)=\wtN_{W,\ell}(G),
\]
and by Lemma~\ref{lem:supported-character-counts}(a) there are \(t^{c(F)-1}\) such
gains.  Their total contribution to \eqref{eq:ie-trace} is therefore
\[
   \sum_{F\subseteq E}
      (-1)^{|E|-|F|}
      t^{-|E|+|F|}\,t^{c(F)-1}\,\wtN_{W,\ell}(G).
\]
Since \(g=|E|-|V|+1\), this equals
\[
   \frac{\wtN_{W,\ell}(G)}{t^g}
   \sum_{F\subseteq E}
      (-1)^{|E|-|F|}
      t^{\beta(F)}
   =
   \frac{\CG(t)}{t^g}\,\wtN_{W,\ell}(G).
\]

Next suppose that \(G\) is nonbipartite and \(t\) is even.  Then the
\(\theta_t\)-class is distinct from the trivial class.  For gains in the
\(\theta_t\)-class,
Proposition~\ref{prop:spectral-input}(b) gives
\[
   \Tr(A_{\gamma,W}^\ell)=(-1)^\ell\wtN_{W,\ell}(G).
\]
By Lemma~\ref{lem:supported-character-counts}(b), such gains in
\(\mathcal G_t(E\setminus F)\) exist exactly when \(G[F]\) is bipartite, and
then there are \(t^{c(F)-1}\) of them.  Hence the \(\theta_t\)-class contribution is
\[
   (-1)^\ell
   \sum_{\substack{F\subseteq E\\ G[F]\ \mathrm{bipartite}}}
      (-1)^{|E|-|F|}
      t^{-|E|+|F|}\,t^{c(F)-1}\,\wtN_{W,\ell}(G)
   =
   (-1)^\ell\frac{\DG(t)}{t^g}\,\wtN_{W,\ell}(G).
\]
In all other cases there is no distinct canonical leading class: either \(t\) is
odd, or \(G\) is bipartite and \(\theta\) is the trivial character.

It remains to bound all other gains.  For each \(F\), the set
\(\mathcal G_t(E\setminus F)\) is finite, and there are only finitely many subsets
\(F\).  By Proposition~\ref{prop:spectral-input}(d), every gain not in the leading
class or classes has spectral radius strictly smaller than \(\rho_{A,W}\).  Because
the collection is finite, there is a common \(\sigma'_A<\rho_{A,W}\) for which
their total contribution is \(O((\sigma'_A)^\ell)\).  Therefore
\[
   \wtN_{W,\ell}^{t\text{-}\NZ}(G)
      =
      \frac{\Phi^A_{G,t}(\ell)}{t^g}\,\wtN_{W,\ell}(G)
         +O((\sigma'_A)^\ell).
\]
Finally, \(A_W\) is an irreducible nonnegative matrix whose period is one when
\(G\) is nonbipartite and two when \(G\) is bipartite.
Perron--Frobenius theory therefore gives
\[
   \wtN_{W,\ell}(G)=
   \begin{cases}
   \rho_{A,W}^\ell+O((\sigma''_A)^\ell),
      & G\ \text{is nonbipartite},\\
   2\rho_{A,W}^\ell+O((\sigma''_A)^\ell),
      & G\ \text{is bipartite and }2\mid \ell,\\
   0, & G\ \text{is bipartite and }2\nmid \ell.
   \end{cases}
\]
for some \(\sigma''_A<\rho_{A,W}\).  Taking
\(\sigma_A=\max\{\sigma'_A,\sigma''_A\}\) gives the two estimates in part~(a).
If \(G\) is bipartite and \(2\nmid\ell\), there are no closed walks of
length \(\ell\), so \(\wtN_{W,\ell}^{t\text{-}\NZ}(G)=0\).
This proves part (a).

For part (b), assume \(g\ge2\).  The nonbacktracking proof is parallel, but we spell
out the substitutions.  Using the twisted weighted nonbacktracking matrix
\(B_{\gamma,W}\), Proposition~\ref{prop:edge-trace-formula} gives, for every
\(S\subseteq E\) and \(\alpha\in M_t(S)\),
\[
   N_{W,t,S}(\alpha,\ell)
      =
      \frac{1}{t^{|S|}}
      \sum_{\gamma\in\mathcal G_t(S)}
         \gamma(\alpha)^{-1}
         \Tr(B_{\gamma,W}^\ell).
\]
Using the same inclusion-exclusion identity as above and putting \(S=E\setminus F\),
we get
\[
   N_{W,\ell}^{t\text{-}\NZ}(G)
      =
      \sum_{F\subseteq E}(-1)^{|E|-|F|}
      \frac{1}{t^{|E|-|F|}}
      \sum_{\gamma\in\mathcal G_t(E\setminus F)}\Tr(B_{\gamma,W}^\ell).
\]
The supported trivial gains contribute, exactly as above,
\[
   \frac{\CG(t)}{t^g}N_{W,\ell}(G).
\]
Now assume that \(t\) is even and that the period character exists.
Lemma~\ref{lem:supported-character-counts}(b) shows that gains in the
\(\eta_t\)-class belonging to \(\mathcal G_t(E\setminus F)\) exist exactly when
\(\eta|_{H_1(G[F],\bbZ)}=1\), and then there are \(t^{c(F)-1}\) of them.
Proposition~\ref{prop:spectral-input}(c) gives
\[
   \Tr(B_{\gamma,W}^\ell)
      =\ee\!\left(\frac{\ell}{2\nu_G}\right)N_{W,\ell}(G).
\]
Their total contribution is therefore
\[
   \frac{\ee(\ell/(2\nu_G))N_{W,\ell}(G)}{t^g}
   \sum_{\substack{F\subseteq E\\
          \eta|_{H_1(G[F],\bbZ)}=1}}
      (-1)^{|E|-|F|}t^{\beta(F)}
   =
   \ee\!\left(\frac{\ell}{2\nu_G}\right)
   \frac{\Deta(t)}{t^g}N_{W,\ell}(G).
\]
When \(\nu_G\mid\ell\), the phase is \((-1)^{\ell/\nu_G}\).  If \(t\) is odd or
the period character does not exist, there is no distinct period class in the
finite dual and hence no second leading contribution.

By Proposition~\ref{prop:spectral-input}(e), every remaining gain has spectral
radius strictly smaller than \(\rho_{B,W}\).  The finite collection of gains in
the inclusion-exclusion sum therefore has a common exponential bound
\(O((\sigma'_B)^\ell)\) for some \(\sigma'_B<\rho_{B,W}\).  Along
\(\nu_G\mid\ell\), we have proved
\[
   N_{W,\ell}^{t\text{-}\NZ}(G)
      =
      \frac{\Phi^B_{G,t}(\ell)}{t^g}N_{W,\ell}(G)+O((\sigma'_B)^\ell).
\]
Finally, the positive weights do not change the support or the period of the
nonbacktracking edge shift.  Perron-Frobenius periodicity on its cyclic core gives
\[
   N_{W,\ell}(G)=
   \begin{cases}
   \nu_G\rho_{B,W}^\ell+O((\sigma''_B)^\ell), & \nu_G\mid \ell,\\
   0, & \nu_G\nmid \ell.
   \end{cases}
\]
for some \(\sigma''_B<\rho_{B,W}\).  Taking
\(\sigma_B=\max\{\sigma'_B,\sigma''_B\}\) gives the estimates in part~(b).
If \(\nu_G\nmid \ell\), there are no circuits of length \(\ell\), so
\(N_{W,\ell}^{t\text{-}\NZ}(G)=0\).  This proves part (b) for every connected
graph of genus at least two, every modulus, and every positive directed weight.
\end{proof}

For connected nonbipartite \(G\), positive directed \(W\), and even \(t\), the
closed-walk part of Theorem~\ref{thm:main}, divided by the total weighted count,
gives
\[
   \lim_{\substack{\ell\to\infty\\2\mid\ell}}
   \frac{\wtN_{W,\ell}^{t\text{-}\NZ}(G)}{\wtN_{W,\ell}(G)}
      =\frac{\CG(t)+\DG(t)}{t^g},
   \qquad
   \lim_{\substack{\ell\to\infty\\2\nmid\ell}}
   \frac{\wtN_{W,\ell}^{t\text{-}\NZ}(G)}{\wtN_{W,\ell}(G)}
      =\frac{\CG(t)-\DG(t)}{t^g}.
\]
Thus, if \(r_{\mathrm{even}}\) and \(r_{\mathrm{odd}}\) denote these two limiting
densities, then
\[
   \frac{\CG(t)}{t^g}
      =\frac{r_{\mathrm{even}}+r_{\mathrm{odd}}}{2},
   \qquad
   \frac{\DG(t)}{t^g}
      =\frac{r_{\mathrm{even}}-r_{\mathrm{odd}}}{2}.
\]
The ordinary flow polynomial supplies the average nowhere-zero density, whereas
the modulating flow polynomial supplies the parity difference.  Thus the
finite-flow decomposition in Proposition~\ref{prop:parity-bias} appears directly
in the asymptotic walk-counting densities.

The preceding proof treats \(t\) as fixed.  For integral abelianization, however, a
single fixed modulus is not enough: a walk whose edge coordinates are nonzero modulo
\(t\) is certainly integrally nowhere-zero, but the converse is not true for any one
choice of \(t\).  The point is that these mod-\(t\) lower bounds improve as
\(t\to\infty\).  The next subsection makes this limiting passage precise.

\subsection{Integral nowhere-zero abelianizations}

The fixed-modulus theorem also gives an approximation result for the integral
nowhere-zero condition.  This condition is stronger in its definition but easier in
the limit: a coefficient which is nonzero modulo \(t\) is certainly nonzero over
\(\bbZ\).  Thus mod-\(t\) nowhere-zero walks give lower approximations to integrally
nowhere-zero walks.  Letting \(t\) grow, the degree estimates for \(\CG(t)\) and
\(\DG(t)\) show that these lower approximations become asymptotically sharp.

\begin{lemma}\label{lem:Phi-large-t}
Let \(G\) be connected and bridgeless.  Then, uniformly in \(\ell\),
\[
   \frac{\Phi^A_{G,t}(\ell)}{t^g}=1+O_G(t^{-1})
   \qquad (t\to\infty).
\]
If \(g\ge2\), then uniformly along \(\nu_G\mid\ell\),
\[
   \frac{\Phi^B_{G,t}(\ell)}{t^g}=1+O_G(t^{-1})
   \qquad (t\to\infty).
\]
\end{lemma}

\begin{proof}
Since \(G\) is connected and bridgeless, the ordinary flow polynomial has degree
\(g\) and leading coefficient \(1\).  Hence
\[
   \frac{\CG(t)}{t^g}=1+O_G(t^{-1}).
\]
If the modulating term does not appear in \(\Phi^A_{G,t}(\ell)\), this proves the
first claim.  If it does appear, then \(G\) is nonbipartite and \(t\) is even.
By Proposition~\ref{prop:canonical-coefficient}, \(\deg\DG\le g-1\), and therefore
\[
   \frac{\DG(t)}{t^g}=O_G(t^{-1}).
\]
This proves the assertion for \(\Phi^A\).  For \(\Phi^B\), the only possible
additional term is \(\Deta(t)\).  Proposition~\ref{prop:period-coefficient} gives
\(\deg\Deta\le g-1\), so the same argument proves the second estimate.
\end{proof}

\begin{corollary}[Integral nowhere-zero approximation]\label{cor:integral-approx}
Let \(G\) be a connected graph and let
\(W:\mathbf E(G)\to\bbR_{>0}\) be a positive directed weight.
\begin{enumerate}[label=(\alph*)]
\item If \(G\) has a bridge, then
\[
   \wtN_{W,\ell}^{\bbZ\text{-}\NZ}(G)=0
   \qquad\text{for every }\ell.
\]
If \(G\) is bridgeless, then, along lengths with \(\wtN_{W,\ell}(G)>0\),
\[
   \frac{\wtN_{W,\ell}^{\bbZ\text{-}\NZ}(G)}{\wtN_{W,\ell}(G)}
      \longrightarrow 1.
\]
\item If \(g\ge2\) and \(G\) has a bridge, then
\[
   N_{W,\ell}^{\bbZ\text{-}\NZ}(G)=0
   \qquad\text{for every }\ell.
\]
If \(g\ge2\) and \(G\) is bridgeless, then, along lengths with
\(N_{W,\ell}(G)>0\)
(equivalently, in the asymptotic range, along \(\nu_G\mid\ell\)),
\[
   \frac{N_{W,\ell}^{\bbZ\text{-}\NZ}(G)}{N_{W,\ell}(G)}
      \longrightarrow 1.
\]
\end{enumerate}
\end{corollary}

\begin{proof}
If \(e\) is a bridge, every closed walk crosses \(e\) equally often in the two
directions.  Hence the \(e\)-coefficient of the integral abelianization is always
zero.  This proves the vanishing assertions for closed walks and for circuits.

Assume now that \(G\) is bridgeless.  For every fixed \(t\), a closed walk which is
nowhere-zero modulo \(t\) has integrally nowhere-zero abelianization.  Therefore
\[
   \wtN_{W,\ell}^{\bbZ\text{-}\NZ}(G)
      \ge \wtN_{W,\ell}^{t\text{-}\NZ}(G).
\]
Along lengths with \(\wtN_{W,\ell}(G)>0\), Theorem~\ref{thm:main}(a) gives
\[
      \frac{\wtN_{W,\ell}^{t\text{-}\NZ}(G)}{\wtN_{W,\ell}(G)}
      =
      \frac{\Phi^A_{G,t}(\ell)}{t^g}+o(1)
      \qquad(\ell\to\infty),
\]
where \(t\) is fixed.  By Lemma~\ref{lem:Phi-large-t}, the right-hand side has
liminf at least \(1-O_G(t^{-1})\).  Since \(t\) can be chosen arbitrarily large and
\(\wtN_{W,\ell}^{\bbZ\text{-}\NZ}(G)\le\wtN_{W,\ell}(G)\), the closed-walk ratio
tends to \(1\).

The proof for circuits is identical.  For fixed \(t\),
\[
   N_{W,\ell}^{\bbZ\text{-}\NZ}(G)
      \ge N_{W,\ell}^{t\text{-}\NZ}(G),
\]
and Theorem~\ref{thm:main}(b), applied along the lengths with
\(N_{W,\ell}(G)>0\), gives
\[
      \frac{N_{W,\ell}^{t\text{-}\NZ}(G)}{N_{W,\ell}(G)}
      =
      \frac{\Phi^B_{G,t}(\ell)}{t^g}+o(1).
\]
Lemma~\ref{lem:Phi-large-t} again lets \(t\to\infty\), and the ratio is bounded
above by \(1\).  This proves the nonbacktracking assertion.
\end{proof}

\subsection{Weighted postman circuits}
\label{sec:weighted-postman}

A \emph{postman circuit} is a circuit that
traverses every unoriented edge of \(G\) at least once.  Define its weighted count
by
\[
   \pc_{W,\ell}(G)
      =
      \sum_{\substack{P\ \mathrm{circuit},\ |P|=\ell\\
                      \text{every }e\in E(G)\text{ occurs in }P}}
         W(P).
\]
The name refers to the Chinese postman problem~\cite{EdmondsJohnson1973}: such a
circuit covers every edge, viewed as a street that must be served.  We impose no
minimal-length condition, so a postman circuit here need not be an optimal postman
tour.
Thus the rooting and orientation conventions are the same as for
\(N_{W,\ell}(G)=\Tr(B_W^{\ell})\).  The corresponding unweighted objects were
called \emph{edge covering cycles} by da Costa and
Policarpo~\cite{DaCostaPolicarpo2018}; their trace count retains the possible
initial oriented edges, as does \(\pc_{\ell}(G)\), before passage to equivalence
classes under cyclic permutation.

Full edge support is weaker than integral nowhere-zero abelianization.  Indeed, a
postman circuit may traverse an edge in both directions with zero net coefficient.
Nevertheless, Corollary~\ref{cor:integral-approx} immediately shows that postman
circuits are asymptotically generic on a bridgeless graph.

\begin{corollary}[Weighted postman-circuit asymptotics]
\label{cor:postman-asymptotic}
Let \(G\) be connected and bridgeless with genus \(g\ge2\), and let
\(W:\mathbf E(G)\to\bbR_{>0}\) be a positive directed weight.  Along lengths for
which \(N_{W,\ell}(G)>0\),
\[
   \frac{\pc_{W,\ell}(G)}{N_{W,\ell}(G)}\longrightarrow1.
\]
Equivalently,
\[
   \pc_{W,\ell}(G)
      =
      \begin{cases}
      \nu_G\rho_{B,W}^{\ell}(1+o(1)),&\nu_G\mid\ell,\\
      0,&\nu_G\nmid\ell.
      \end{cases}
\]
\end{corollary}

\begin{proof}
If \(P^{\mathrm{ab}}\) is nowhere zero over \(\bbZ\), then every edge occurs in
\(P\).  Positivity of the directed weight therefore gives
\[
   N_{W,\ell}^{\bbZ\text{-}\NZ}(G)
      \le \pc_{W,\ell}(G)
      \le N_{W,\ell}(G).
\]
Divide by \(N_{W,\ell}(G)\) and apply
Corollary~\ref{cor:integral-approx}(b).  The second formula follows from
Perron--Frobenius periodicity for \(B_W\).
\end{proof}

\ref{app:postman-pf} gives a direct Perron--Frobenius proof based on the
exact weighted inclusion-exclusion formula for postman circuits.  That argument
also shows that bridgelessness can be weakened to the absence of degree-one
vertices.  The underlying unweighted inclusion-exclusion formula first appeared in
da Costa and Policarpo~\cite[Theorem~2.1]{DaCostaPolicarpo2018}.  A directed-support
version of the same transform was later obtained independently in the periodic-orbit algorithm of
Band, Exner, Goel, and Strauss~\cite{BandExnerGoelStrauss2026}, where it is used in
particular to count Eulerian circuits.

\section{Tutte and matroidal structure of the character family}
\label{sec:matroidal-structure}

The Fourier description produces the character-indexed family \(D_{G,\zeta}\).
Its canonical and period members are the two distinguished modulating flow
polynomials selected by the asymptotic theorem.
We first place the entire family in the gain-semimatroid Tutte framework and derive
the structural properties common to all its members.  We then specialize to the
canonical character, whose parity condition makes \(D_G\) an invariant of the cycle
matroid.  Finally, characters of order at most two, including the period character
when it exists, are placed in the signed-Tutte framework, which also supplies the
natural deletion-contraction rule.

For comparison, let \(M(G)\) be the cycle matroid of \(G\), let \(M(G)^*\) be its
dual, and write
\[
   \chi_N(t)=\sum_{A\subseteq E}(-1)^{|A|}
      t^{r_N(E)-r_N(A)}
\]
for the characteristic polynomial of a matroid \(N\) on \(E\).  The classical
cographic formulation of the flow polynomial is
\begin{equation}\label{eq:classical-cographic}
   \CG(t)=\chi_{M(G)^*}(t).
\end{equation}
It is an immediate consequence of Tutte's flow specialization and matroid
duality~\cite{Tutte1954}.  This identity is background for the different
cycle-matroid formulation of \(D_G\) below.

\subsection{Semimatroid and cycle-matroid structure}

Fix a full complex character
\[
   \zeta:H_1(G,\bbZ)\longrightarrow\bbC^\ast,
\]
and choose a \(\bbC^\ast\)-gain graph \(\Phi_\zeta\) representing the switching
class determined by \(\zeta\); equivalently, the gain around every oriented cycle
\(C\) is \(\zeta([C])\).  By
Lemma~\ref{lem:balanced-dichromatic-specialization}, \(D_{G,\zeta}\) is a
specialization of the balanced dichromatic polynomial of \(\Phi_\zeta\).  We now
realize its balanced edge sets as the central sets of a semimatroid.  Recall the
defining expansion
\[
   D_{G,\zeta}(q)
      =
      \sum_{\substack{F\subseteq E\\
         \zeta|_{H_1(G[F],\bbZ)}=1}}
      (-1)^{|E|-|F|}q^{\beta(F)}.
\]
Put
\[
   \mathcal B_{G,\zeta}
      =\{F\subseteq E:\zeta|_{H_1(G[F],\bbZ)}=1\}.
\]
These are exactly the balanced edge sets of the gain graph representing \(\zeta\).
Balanced edge sets of a biased graph form the central sets of a semimatroid; see
Ardila~\cite{Ardila2007Semimatroids} and Fu~\cite[Proposition~3.5]{Fu2026Semimatroids}.
Write
\[
   \mathcal S_{G,\zeta}
      =(E,\mathcal B_{G,\zeta},r_{\mathcal S}),
   \qquad
   r_{\mathcal S}(F)=|V|-c(F).
\]
Every spanning forest is balanced, so
\(r(\mathcal S_{G,\zeta})=|V|-c(E)=r_{M(G)}(E)\).  The semimatroid Tutte polynomial is
\begin{equation}\label{eq:semimatroid-tutte}
   T_{\mathcal S_{G,\zeta}}(x,y)=
      \sum_{F\in\mathcal B_{G,\zeta}}
      (x-1)^{r(\mathcal S_{G,\zeta})-r_{\mathcal S}(F)}
      (y-1)^{|F|-r_{\mathcal S}(F)}.
\end{equation}

\begin{proposition}[Structure of the character-indexed family]
\label{prop:semimatroid-specialization}
For every finite graph \(G\) and every \(\zeta\in\mathcal X_{\bbC}(G)\),
the following hold.
\begin{enumerate}[label=(\alph*)]
\item The \(\zeta\)-member is the gain-semimatroid specialization
   \begin{equation}\label{eq:semimatroid-specialization}
      D_{G,\zeta}(q)
         =(-1)^{\beta(G)}T_{\mathcal S_{G,\zeta}}(0,1-q).
   \end{equation}

\item Its constant term is
   \begin{equation}\label{eq:constant-term}
      D_{G,\zeta}(0)=\CG(0)=(-1)^{\beta(G)}T_G(0,1),
   \end{equation}
   and
   \[
      D_{G,\zeta}(q)\not\equiv0
         \quad\Longleftrightarrow\quad
      G\ \text{is bridgeless}.
   \]

\item If \(G\) is bridgeless and \(d=\deg D_{G,\zeta}\), then, writing
   \([q^k]f(q)\) for the coefficient of \(q^k\) in \(f(q)\),
   \begin{equation}\label{eq:strict-coefficient-signs}
      (-1)^{\beta(G)+k}[q^k]D_{G,\zeta}(q)>0
         \qquad(0\le k\le d).
   \end{equation}
   Thus the coefficients strictly alternate in sign and have no internal zeros.
\end{enumerate}
\end{proposition}

\begin{proof}
For part~(a), put \(r=r(\mathcal S_{G,\zeta})\) and \(m=|E|\), so that
\(\beta(G)=m-r\).  For \(F\in\mathcal B_{G,\zeta}\), its contribution to
\(T_{\mathcal S_{G,\zeta}}(0,1-q)\) is
\[
   (-1)^{r-r_{\mathcal S}(F)}(-q)^{|F|-r_{\mathcal S}(F)}
      =(-1)^{r+|F|}q^{\beta(F)}.
\]
Multiplication by \((-1)^{\beta(G)}=(-1)^{m-r}\) changes the sign to
\((-1)^{m+|F|}=(-1)^{m-|F|}\).  Summing over the balanced sets gives
\eqref{eq:general-character-polynomial}.

For part~(b), the constant term in \eqref{eq:general-character-polynomial} is
obtained from the sets \(F\) satisfying \(\beta(F)=0\).  These are precisely the
forests, and every character is trivial on their zero homology groups.  Hence
\(D_{G,\zeta}(0)=\CG(0)\), and the last equality in \eqref{eq:constant-term} follows
from the usual subset expansion of \(T_G\).

If \(e\) is a bridge, pair every balanced set \(F\) not containing \(e\) with
\(F\cup\{e\}\).  Adding \(e\) changes neither the cycle rank nor the restricted
character, while it reverses the sign in
\eqref{eq:general-character-polynomial}; hence all terms cancel and
\(D_{G,\zeta}=0\).  Conversely, suppose that \(G\) is bridgeless.  The loops of the
cographic matroid \(N=M(G)^*\) are precisely the bridges of \(G\), so \(N\) is
loopless and has rank \(r_N(E)=\beta(G)\).  Let \(L(N)\) be the lattice of
flats of \(N\), and let \(\mu_N\) denote its M\"obius function.  Since \(N\) is
loopless, the minimum and maximum elements of \(L(N)\) are \(\varnothing\) and
\(E\).  Thus
\[
   \chi_N(t)=\sum_{F\in L(N)}
      \mu_N(\varnothing,F)t^{r_N(E)-r_N(F)},
\]
and hence
\[
   \chi_N(0)=\mu_N(\varnothing,E).
\]
The lattice \(L(N)\) is geometric, and the sign theorem for geometric lattices
gives
\[
   (-1)^{\beta(G)}\mu_N(\varnothing,E)>0
\]
(see Rota~\cite[Section~7]{Rota1964Mobius}).  Thus
\(\chi_{M(G)^*}(0)\ne0\), and the classical cographic identity
\eqref{eq:classical-cographic} together with \eqref{eq:constant-term} gives
\(D_{G,\zeta}(0)\ne0\).

For part~(c), the semimatroid activity expansion has nonnegative integral
coefficients~\cite{Ardila2007Semimatroids}.  Write
\[
   T_{\mathcal S_{G,\zeta}}(0,y)=\sum_{j=0}^{s}a_jy^j,
   \qquad a_j\ge0,\quad a_s>0.
\]
Part~(a) shows that
\(d=\deg D_{G,\zeta}=s\) and gives
\[
   D_{G,\zeta}(q)
      =\sum_{k=0}^{d}(-1)^{\beta(G)+k}
         \left(\sum_{j=k}^{d}\binom{j}{k}a_j\right)q^k.
\]
For every \(k\le d\), the inner sum is positive because it contains
\(\binom{d}{k}a_d\).
\end{proof}

Taking \(\zeta=\theta\), and taking \(\zeta=\eta\) whenever the componentwise
period character exists, gives the same bridgelessness criterion and coefficient
conclusions for \(\DG\) and \(\Deta\), respectively.

The nonvanishing assertion in
Proposition~\ref{prop:semimatroid-specialization} concerns the polynomial, not
each prescribed integer.  Thus bridgelessness guarantees
\(\DG(t)\not\equiv0\), but a particular value \(\DG(t_0)\) may vanish through
cancellation even when \(\CG(t_0)>0\).  Example~\ref{ex:evaluation-zeros} exhibits
this phenomenon for the Petersen graph and for \(K_5\) with two adjacent edges
deleted.  Such zeros record balance between the two canonical classes at \(t_0\);
they are not flow-existence obstructions.

\begin{remark}[Cycle-hypertorus stratification]
Although \(\mathcal X_{\bbC}(G)\) is generally a positive-dimensional torus, a
fixed graph produces only finitely many polynomials \(D_{G,\zeta}\).  Indeed, the
polynomial records only which cycles \(C\) satisfy \(\zeta([C])=1\).
Thus the character torus is stratified by the finitely many cycle hypertori,
one for each simple cycle \(C\),
\[
   \mathcal H_C
      =\{\zeta\in\mathcal X_{\bbC}(G):\zeta([C])=1\}.
\]
By a \emph{Zariski-generic character} we mean a character in the Zariski-open
dense complement of \(\bigcup_C\mathcal H_C\); equivalently,
\(\zeta([C])\ne1\) for every simple cycle \(C\).  Its balanced edge sets are
therefore precisely the forests, and
\[
   D_{G,\zeta}(q)=\CG(0).
\]
For example, on a cycle,
\[
   D_{C_n,\mathbf 1}(q)=q-1,
   \qquad
   D_{C_n,\zeta}(q)=-1\quad(\zeta\ne\mathbf 1).
\]
For the theta multigraph consisting of three parallel edges, the three possibilities,
according as all, exactly one, or none of its cycles are balanced, are
\[
   (q-1)(q-2),\qquad 2-q,\qquad 2.
\]
Two balanced cycles force the third to be balanced.  The full character family is
therefore finite at the polynomial level.  The spectral results in
Sections~\ref{sec:trace-selection} and~\ref{sec:main-proof}
explain why only its trivial and canonical members enter the leading vertex-walk
terms and why, in genus at least two, only its trivial and possible period members
enter the leading nonbacktracking terms.  Genus one is the exceptional edge case.
\end{remark}

The canonical character admits a further formulation involving only the cycle
matroid.  Put \(M=M(G)\), with rank function
\[
   r_M(F)=|V|-c(F),\qquad F\subseteq E,
\]
so that \(\beta(F)=|F|-r_M(F)\).  Call a restriction \(N|F\) of a matroid
\(N\) \emph{even} if every circuit of \(N|F\) has even cardinality, and define
\begin{equation}\label{eq:matroidal-D}
   D_N(t)=
      \sum_{\substack{F\subseteq E\\N|F\ \mathrm{even}}}
      (-1)^{|E|-|F|}t^{|F|-r_N(F)}.
\end{equation}
Unlike the classical cographic identity \eqref{eq:classical-cographic}, the
canonical polynomial is obtained from restrictions of the cycle matroid \(M(G)\).

\begin{lemma}[Canonical cycle-matroid formulation]
\label{prop:matroidal-formulations}
For every finite graph \(G\),
\[
   \DG(t)=D_{M(G)}(t).
\]
In particular, \(\DG(t)\) depends only on the cycle matroid of \(G\).
\end{lemma}

\begin{proof}
The circuits of \(M(G)|F\) are precisely the edge sets of the cycles of \(G[F]\),
with a loop regarded as a cycle of length one.  Hence \(G[F]\) is bipartite if and
only if \(M(G)|F\) is even.  Since \(|F|-r_M(F)=\beta(F)\), the identity follows
from \eqref{eq:matroidal-D} and \eqref{eq:modulating-flow-poly}.
\end{proof}

\subsection{Order-two characters and signed-Tutte structure}

The gain-semimatroid formulation above applies uniformly to every character, but
the two characters distinguished by the asymptotic theorem have additional
order-two structure: the canonical character, and the period character whenever it
exists.  More generally, every character of order at most two admits a
\(\{\pm1\}\)-gain representative and hence a signed-graph model.  This reduction
places the two distinguished polynomials in a common established framework and,
unlike the general complex-gain formulation, supplies a graph-theoretic
deletion-contraction rule.  We therefore develop the signed model for all order-two
characters before specializing to the canonical and period characters.

Let \(\Sigma=(G,\sigma)\) be a signed graph, where
\(\sigma:E(G)\to\{\pm1\}\).  An edge set \(F\subseteq E(G)\) is \emph{balanced}
if every cycle of \(\Sigma[F]\) has positive sign.  Define
\[
   \begin{aligned}
   D_\Sigma(t)
      &=
      \sum_{\substack{F\subseteq E(G)\\ \Sigma[F]\ \mathrm{balanced}}}
      (-1)^{|E(G)|-|F|}t^{|F|-|V(G)|+c(F)}  \\
      &=
      \sum_{\substack{F\subseteq E(G)\\ \Sigma[F]\ \mathrm{balanced}}}
      (-1)^{|E(G)|-|F|}t^{\beta(F)}.
   \end{aligned}
\]
If \(\Sigma_\zeta\) represents a character \(\zeta\) of order at most two, then
the definitions give
\[
   D_{G,\zeta}(t)=D_{\Sigma_\zeta}(t).
\]
Characters of higher order remain in the gain-semimatroid framework and do not, in
general, have a signed-graph representative.

This construction belongs to an established signed-Tutte framework.  The
trivariate signed Tutte polynomial of Goodall--Litjens--Regts--Vena has evaluations
that count nowhere-zero flows on signed graphs~\cite{GoodallLitjensRegtsVena2020}.
They also identify it, after normalization, with the joint Tutte polynomial of the
ordered matroid pair
\[
   \bigl(M(G),F(\Sigma)\bigr),
\]
where \(F(\Sigma)\) is the frame matroid.  The specialization below connects this
framework with the canonical Fourier coefficient and the spectral parity term.

For comparison, the signed flow polynomials of Ren--Qian count nowhere-zero flows
over finite abelian groups with prescribed two-torsion rank~\cite{RenQian2019}.
Their \(d=0\) subset expansion retains both balanced and unbalanced signed
subgraphs and, after reindexing the deletion set, has the form
\[
   F_0(\Sigma,x)
      =
      \sum_{H\subseteq E(G)}(-1)^{|E(G)|-|H|}x^{\beta_\Sigma(H)},
\]
where
\[
   \beta_\Sigma(H)
      =
      \sum_{\omega\in\pi_0(\Sigma[H])}
      \bigl(|E(\omega)|-|V(\omega)|
            +\mathbf 1_{\{\omega\ \mathrm{balanced}\}}\bigr).
\]
For a balanced subgraph \(F\), this exponent is
\(|F|-|V(G)|+c(F)\).  Thus \(D_\Sigma(t)\) is the balanced-subgraph part selected
from the signed-flow setting.  It is generally distinct from \(F_0(\Sigma,x)\),
although the two coincide when \(\Sigma\) is balanced.  The next proposition gives
the signed-Tutte specialization, and Proposition~\ref{prop:signed-deletion-contraction}
below gives its recursive form.

\begin{proposition}[Signed-Tutte specialization]\label{prop:signed-tutte}
Let \(\Sigma=(G,\sigma)\) be a connected signed graph of genus \(g\).  If
\(T_{\Sigma}(X,Y,Z)\) denotes the trivariate Tutte polynomial of a signed graph in the
sense of
Goodall--Litjens--Regts--Vena~\cite{GoodallLitjensRegtsVena2020}, then
\begin{equation}\label{eq:signed-tutte-specialization}
   D_\Sigma(t)=(-1)^g T_{\Sigma}(0,1-t,1).
\end{equation}
In particular, \(\DG(t)=D_{\Sigma_-}(t)\), where \(\Sigma_-\) is the all-negative
signing of \(G\).
\end{proposition}

\begin{proof}
In the subset expansion of \(T_{\Sigma}\), evaluating at \(Z=1\) kills
exactly those subsets which contain an unbalanced connected component, and hence
retains exactly the balanced spanning subgraphs.  For such a subset \(F\), the
surviving signed-Tutte term at \((0,1-t,1)\) is
\[
   (-1)^{c(F)-1}(-t)^{\beta(F)}
      =
   (-1)^g(-1)^{|E|-|F|}t^{\beta(F)},
\]
because \(g=|E|-|V|+1\).  This proves the formula.
For the all-negative signing, a subgraph is balanced exactly when its underlying graph
is bipartite, so \(D_{\Sigma_-}(t)=\DG(t)\).
\end{proof}

\begin{corollary}[Period signed-Tutte specialization]
\label{cor:period-signed-tutte}
Suppose that \(G\) is connected of genus at least two and its period character
\(\eta\) exists.  Let \(\Sigma_\eta=(G,\sigma_\eta)\) be any signing representing
the corresponding class in \(H^1(G,\bbZ/2\bbZ)\).  Then
\[
   \Deta(t)=D_{\Sigma_\eta}(t)
      =(-1)^gT_{\Sigma_\eta}(0,1-t,1).
\]
This polynomial depends only on the switching class representing \(\eta\).
\end{corollary}

\begin{proof}
Choose the signing so that, for every cycle \(C\),
\[
   \eta(C)=\prod_{e\in C}\sigma_\eta(e).
\]
Then \(\eta\) is trivial on \(H_1(G[F],\bbZ)\) exactly when every cycle of
\(\Sigma_\eta[F]\) has positive sign, that is, exactly when
\(\Sigma_\eta[F]\) is balanced.  The defining expansions of \(\Deta\) and
\(D_{\Sigma_\eta}\) therefore agree.  Proposition~\ref{prop:signed-tutte} gives
the specialization.  Changing the cochain representative of \(\eta\) switches the
signing, and switching preserves balanced edge sets.
\end{proof}

The signed-graph formulation is also the natural setting for deletion-contraction.
Contracting an edge in the all-negative signing can produce a mixed signing, so the
class of all-negative signings is not closed under contraction.  The
balanced-subgraph specialization \(D_\Sigma(t)\) is closed under the signed
deletion-contraction operations.  Specializing the deletion-contraction recurrence
of Goodall--Litjens--Regts--Vena~\cite[Theorem~4.1]{GoodallLitjensRegtsVena2020}
gives the following recurrence for \(D_\Sigma\).  As in their framework, a non-loop
edge is switched to positive before contraction, whereas a negative loop is not
contracted.  We include a direct balanced-subgraph proof to make the resulting signs
and loop cases explicit.

\begin{proposition}[Signed deletion-contraction]\label{prop:signed-deletion-contraction}
Let \(\Sigma\) be a signed graph and let \(e\in E(\Sigma)\).
\begin{enumerate}[label=(\alph*)]
\item If \(e\) is not a loop, then
\[
   D_\Sigma(t)=D_{\Sigma/e}(t)-D_{\Sigma\setminus e}(t),
\]
where signed contraction is defined by first switching, if necessary, so that \(e\)
is positive, and then contracting \(e\).
\item If \(e\) is a positive loop, then
\[
   D_\Sigma(t)=(t-1)D_{\Sigma\setminus e}(t).
\]
\item If \(e\) is a negative loop, then
\[
   D_\Sigma(t)=-D_{\Sigma\setminus e}(t).
\]
\end{enumerate}
\end{proposition}

\begin{proof}
Switching a signed graph preserves balanced edge sets, so \(D_\Sigma(t)\) is
switching invariant.  Suppose first that \(e\) is not a loop, and switch so that
\(e\) is positive.  Split the defining expansion for \(D_\Sigma(t)\) according to
whether a balanced set contains \(e\).  The sets not containing \(e\) contribute
\[
   -D_{\Sigma\setminus e}(t),
\]
because the cycle rank \(\beta(F)\) is unchanged while the ambient number of edges
has decreased by one in \(\Sigma\setminus e\).

For the sets containing \(e\), write them as \(F\cup\{e\}\).  Since \(e\) is
positive, balanced sets \(F\cup\{e\}\) in \(\Sigma\) are in bijection with balanced
sets \(F\) in \(\Sigma/e\).  Moreover
\[
   |F\cup\{e\}|-|V(\Sigma)|+c(F\cup\{e\})
      =
   |F|-|V(\Sigma/e)|+c_{\Sigma/e}(F),
\]
and the signs in the two subset expansions agree.  Hence the sets containing \(e\)
contribute \(D_{\Sigma/e}(t)\), proving (a).

If \(e\) is a positive loop, the balanced sets containing \(e\) are exactly
\(F\cup\{e\}\), where \(F\) is balanced in \(\Sigma\setminus e\).  Their exponent is
larger by one than the exponent of \(F\).  Thus the contribution from sets
containing \(e\) is \(tD_{\Sigma\setminus e}(t)\), while the contribution from sets
not containing \(e\) is \(-D_{\Sigma\setminus e}(t)\).  This gives (b).  If \(e\) is
a negative loop, no balanced set contains \(e\), so only the latter contribution
remains, giving (c).
\end{proof}

\begin{remark}[Why signed contraction is necessary]\label{ex:signed-contraction}
For an ordinary graph \(G\) and a non-loop edge \(e\),
\[
   \DG(t)=D_{\Sigma_-/e}(t)-D_{G\setminus e}(t),
\]
where the contraction term is the balanced-subgraph polynomial of a signed graph.  It
is generally not \(D_{G/e}(t)\).  For example, if \(G=C_3\), then \(\DG(t)=-1\) and
\(D_{G\setminus e}(t)=0\), whereas ordinary contraction gives
\[
   D_{G/e}(t)=t-1.
\]
The signed contraction \(\Sigma_-/e\), which has two parallel edges of opposite
signs, gives the correct value \(D_{\Sigma_-/e}(t)=-1\).
The same qualification applies to \(\Deta\):
Proposition~\ref{prop:signed-deletion-contraction} may be applied to
\(\Sigma_\eta\), but the inherited signing on a deletion or contraction need not
represent the intrinsic period character of the minor.
\end{remark}

The results of this section separate the universal structure from the special
features of \(D_G\).  Every \(D_{G,\zeta}\) is a gain-semimatroid Tutte
specialization with the same polynomial nonvanishing criterion and strict coefficient
alternation.  The canonical member is additionally a cycle-matroid invariant and the
all-negative signed-Tutte specialization.  The period member belongs to the
signed-Tutte framework through Corollary~\ref{cor:period-signed-tutte}, but we make no
claim that it is an invariant of the ordinary cycle matroid.  We now turn to examples
that illustrate these structural distinctions together with the parity effects
selected by the walk dynamics.

\section{Examples and comparisons}
\label{sec:examples}

This section records several computations illustrating the roles played by
\(\DG(t)\) and \(\Deta(t)\).  Cycles provide a transparent calibration.  The
\(K_4\) computations then compare direct nonbacktracking enumeration with the
unweighted and weighted asymptotic formulas.  The bipartite graph
\(\Theta_{2,2,4}\) isolates the period correction, and the final comparison shows
that evaluation-specific zeros of \(\DG\) represent cancellation rather than
failure of nowhere-zero flows.  Together the examples exhibit complete parity
concentration, two distinct limiting densities, no leading parity bias, and a
nonzero bias with both parity classes present.

\begin{example}[Cycles]
Let \(C_n\) be the cycle of length \(n\).  Then
\[
   C_{C_n}(t)=t-1,
\]
whereas
\[
   D_{C_n}(t)
      =
      \begin{cases}
      t-1, & n\ \mathrm{even},\\
      -1, & n\ \mathrm{odd}.
      \end{cases}
\]
Indeed, if \(n\) is even then \(C_n\) is bipartite and every spanning subgraph is
bipartite.  If \(n\) is odd, the only non-bipartite spanning subgraph is the full
cycle, so the \(t\)-term in the ordinary flow polynomial is removed and only the
constant contribution from proper subgraphs remains.  This example already shows
that \(D_G(t)\) is not a second copy of the flow polynomial: it separates the
bipartite and non-bipartite parts of the cycle structure.
\end{example}

\begin{example}[\(K_4\) and explicit nonbacktracking counts]\label{ex:k4-counts}
For the complete graph \(K_4\),
the classical flow polynomial and the defining expansion
\eqref{eq:modulating-flow-poly} give
\[
   C_{K_4}(t)=(t-1)(t-2)(t-3),
   \qquad
   D_{K_4}(t)=3(t-2).
\]
The polynomial values needed for the comparison are
\[
\begin{array}{c|ccc}
t & 4 & 5 & 6\\ \hline
C_{K_4}(t) & 6 & 24 & 60\\
D_{K_4}(t) & 6 & 9 & 12
\end{array}
\]
For \(t=4,5,6\), set
\[
   R_t(\ell):=\frac{N_{\ell}^{t\text{-}\NZ}(K_4)}{N_\ell(K_4)} .
\]
The leading densities predicted by Theorem~\ref{thm:main} are therefore the
following subsequential limits:
\[
\begin{array}{c|cc}
t & \ell\ \mathrm{even} & \ell\ \mathrm{odd}\\ \hline
4 & 3/16 & 0\\
5 & 24/125 & 24/125\\
6 & 1/3 & 2/9
\end{array}
\]
Although \(D_{K_4}(5)=9\) as a polynomial value, at the odd modulus \(t=5\) the
canonical correction is absent.  Both parities therefore have the same
leading density \(C_{K_4}(5)/5^3=24/125\).  Thus \(t=4\) gives complete parity
concentration, \(t=6\) gives a genuine two-level parity oscillation, and \(t=5\)
gives a single limiting level.

The following table gives representative exact counts and the corresponding ratios
for small and intermediate lengths; the full data through length \(40\) are recorded
in the accompanying notebook.
\begin{center}
\resizebox{\textwidth}{!}{%
\(
\begin{array}{c|r|rr|rr|rr}
\ell & N_\ell(K_4)
& \multicolumn{2}{c|}{t=4}
& \multicolumn{2}{c|}{t=5}
& \multicolumn{2}{c}{t=6}\\
& & N_{\ell}^{4\text{-}\NZ}(K_4) & R_4(\ell)
& N_{\ell}^{5\text{-}\NZ}(K_4) & R_5(\ell)
& N_{\ell}^{6\text{-}\NZ}(K_4) & R_6(\ell)\\ \hline
14 & 16128 & 2688 & 0.166667 & 3696 & 0.229167 & 3696 & 0.229167\\
15 & 31944 & 0 & 0 & 5040 & 0.157776 & 5040 & 0.157776\\
18 & 259152 & 50112 & 0.193369 & 57888 & 0.223375 & 79488 & 0.306723\\
19 & 526680 & 0 & 0 & 88464 & 0.167965 & 106704 & 0.202597\\
24 & 16759272 & 3129984 & 0.186761 & 3397536 & 0.202726 & 5428224 & 0.323894\\
25 & 33541200 & 0 & 0 & 6129600 & 0.182748 & 7215600 & 0.215126\\
30 & 1073772096 & 201365760 & 0.187531 & 211833120 & 0.197279 & 355126320 & 0.330728\\
31 & 2147211528 & 0 & 0 & 401684112 & 0.187072 & 472258464 & 0.219940\\
36 & 68720897880 & 12886058880 & 0.187513 & 13377555648 & 0.194665 & 22845258144 & 0.332435
\end{array}
\)}
\end{center}
Figure~\ref{fig:k4-ratios} plots the same ratios for \(10\le \ell\le 40\).  Thus
the table records exact numerical values, while the figure displays the approach to
the limiting levels predicted by Theorem~\ref{thm:main}.
The denominator \(N_\ell(K_4)\) is computed as a nonbacktracking trace.  Let \(B\)
be the ordinary Hashimoto matrix of \(K_4\).  The adjacency spectrum of \(K_4\) is
\(3,-1,-1,-1\), and the Ihara--Bass spectral relation~\cite{Terras2011} gives
\[
   \operatorname{spec}(B)
      =
      \{2,1,1,1,-1,-1,\alpha,\alpha,\alpha,
        \overline{\alpha},\overline{\alpha},\overline{\alpha}\},
   \qquad
   \alpha=\frac{-1+i\sqrt7}{2}.
\]
Consequently
\[
   N_\ell(K_4)=\Tr(B^\ell)
      =
      2^\ell+3+2(-1)^\ell
      +3\bigl(\alpha^\ell+\overline{\alpha}^{\,\ell}\bigr).
\]
For example, this gives \(N_{14}(K_4)=16128\) and
\(N_{36}(K_4)=68720897880\), the first and last denominators in the displayed
table.  The nowhere-zero numerators are obtained by dynamic programming on the
oriented-edge state space: the state records the current oriented edge and the
residue vector in \((\bbZ/t\bbZ)^{E(K_4)}\), and the final step enforces the cyclic
nonbacktracking condition.  This is an explicit enumeration of rooted, oriented
circuits, not an evaluation of the limiting formula.

The exact enumeration is already close to the theoretical limits.  Comparing the
last even and odd lengths used in Figure~\ref{fig:k4-ratios} gives
\[
\begin{array}{c|cc|cc}
t & \multicolumn{2}{c|}{\ell\ \mathrm{even}}
  & \multicolumn{2}{c}{\ell\ \mathrm{odd}}\\
& \mathrm{theory} & \mathrm{counting}\ \ell=40
& \mathrm{theory} & \mathrm{counting}\ \ell=39\\ \hline
4 & 0.187500 & 0.187501 & 0 & 0\\
5 & 0.192000 & 0.193655 & 0.192000 & 0.190134\\
6 & 0.333333 & 0.332898 & 0.222222 & 0.221697
\end{array}
\]
Thus the persistent parity splitting at \(t=4\) and \(t=6\) is visible both in the
theorem and in the exact counts.  At \(t=5\), by contrast, the even and odd
subsequences approach the same limit, so the remaining oscillation in the displayed
finite lengths is a lower-order effect.

\begin{figure}[H]
\centering
\begin{tikzpicture}
\begin{axis}[
   width=0.86\textwidth,
   height=0.46\textwidth,
   xlabel={length \(\ell\)},
   ylabel={\(R_t(\ell)\)},
   xmin=10, xmax=40,
   ymin=0, ymax=0.42,
   xtick={10,15,20,25,30,35,40},
   ytick={0,0.1,0.2,0.3,0.4},
   grid=both,
   legend style={at={(0.5,-0.24)},anchor=north,legend columns=3},
]
\addplot[black!35, densely dashed, mark=none, forget plot] coordinates {
(10,0.187500) (40,0.187500)};
\addplot[black!35, densely dashed, mark=none, forget plot] coordinates {
(10,0.192000) (40,0.192000)};
\addplot[black!35, densely dashed, mark=none, forget plot] coordinates {
(10,0.222222) (40,0.222222)};
\addplot[black!35, densely dashed, mark=none, forget plot] coordinates {
(10,0.333333) (40,0.333333)};
\addplot+[mark=*, mark size=1.3pt] coordinates {
(10,0.400000) (11,0.000000) (12,0.218182) (13,0.000000)
(14,0.166667) (15,0.000000) (16,0.177969) (17,0.000000)
(18,0.193369) (19,0.000000)
(20,0.187040) (21,0.000000) (22,0.186257) (23,0.000000)
(24,0.186761) (25,0.000000) (26,0.187670) (27,0.000000)
(28,0.187520) (29,0.000000) (30,0.187531) (31,0.000000)
(32,0.187537) (33,0.000000) (34,0.187486) (35,0.000000)
(36,0.187513) (37,0.000000) (38,0.187499) (39,0.000000)
(40,0.187501)};
\addlegendentry{\(t=4\)}
\addplot+[mark=square*, mark size=1.3pt] coordinates {
(10,0.400000) (11,0.000000) (12,0.218182) (13,0.214286)
(14,0.229167) (15,0.157776) (16,0.218156) (17,0.173913)
(18,0.223375) (19,0.167965)
(20,0.207568) (21,0.175660) (22,0.205764) (23,0.180361)
(24,0.202726) (25,0.182748) (26,0.201375) (27,0.184605)
(28,0.199065) (29,0.185901) (30,0.197279) (31,0.187072)
(32,0.196346) (33,0.188220) (34,0.195399) (35,0.188986)
(36,0.194665) (37,0.189646) (38,0.194113) (39,0.190134)
(40,0.193655)};
\addlegendentry{\(t=5\)}
\addplot+[mark=triangle*, mark size=1.6pt] coordinates {
(10,0.400000) (11,0.000000) (12,0.218182) (13,0.214286)
(14,0.229167) (15,0.157776) (16,0.298529) (17,0.211180)
(18,0.306723) (19,0.202597)
(20,0.324354) (21,0.209440) (22,0.321797) (23,0.210312)
(24,0.323894) (25,0.215126) (26,0.327340) (27,0.217450)
(28,0.329647) (29,0.218959) (30,0.330728) (31,0.219940)
(32,0.331420) (33,0.220684) (34,0.331978) (35,0.221181)
(36,0.332435) (37,0.221451) (38,0.332710) (39,0.221697)
(40,0.332898)};
\addlegendentry{\(t=6\)}
\end{axis}
\end{tikzpicture}
\caption{Exact ratios \(R_t(\ell)\) for rooted, oriented circuits
in \(K_4\).  The evaluations \(t=4\) and \(t=6\) exhibit parity oscillations,
whereas the odd evaluation \(t=5\) has a single limiting density.  The dashed
horizontal lines mark the limiting levels predicted by Theorem~\ref{thm:main}.}
\label{fig:k4-ratios}
\end{figure}
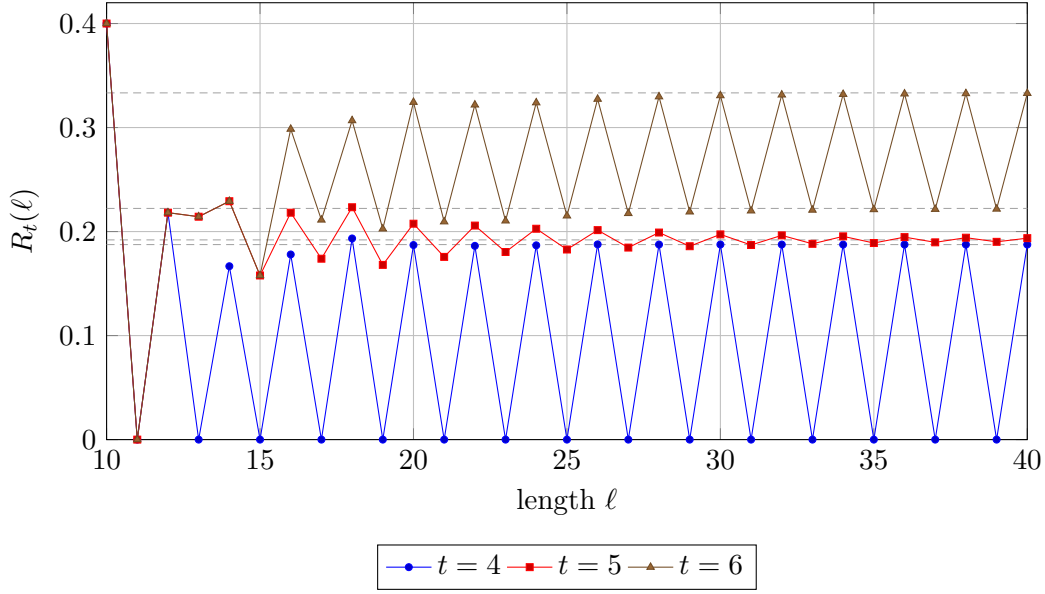

\ref{app:dynamic-count} records the dynamic enumeration used for the
exact numerators in the table.  The accompanying notebook implements the same
recurrence and generates the displayed data and plot.
\end{example}

\begin{example}[An asymmetric directed weight on \(K_4\)]
\label{ex:weighted-k4}
Label the vertices of \(K_4\) by \(0,1,2,3\), and assign the positive directed
weight
\[
   W(0\!\to\!1)=2,
   \qquad
   W(\mathbf e)=1
   \quad\text{for every other }\mathbf e\in\mathbf E(K_4).
\]
In particular, \(W(0\!\to\!1)\ne W(1\!\to\!0)\), so this weight is not induced
by an undirected edge weighting.  If \(m_{01}(P)\) denotes the number of
traversals of \(0\!\to\!1\) by a circuit \(P\), then
\[
   W(P)=2^{m_{01}(P)}.
\]
Thus the weighted sums distinguish circuits according to their directed use of a
single edge.

Consider the even modulus \(t=6\).  Since \(K_4\) is nonbipartite and
has genus \(g=3\), Theorem~\ref{thm:main}(b), together with
\(C_{K_4}(6)=60\) and \(D_{K_4}(6)=12\), gives
\[
   \frac{N_{W,\ell}^{6\text{-}\NZ}(K_4)}{N_{W,\ell}(K_4)}
      \longrightarrow
      \begin{cases}
         \dfrac{60+12}{6^3}=\dfrac13,& \ell\ \mathrm{even},\\[1ex]
         \dfrac{60-12}{6^3}=\dfrac29,& \ell\ \mathrm{odd}.
      \end{cases}
\]
The coefficients in these limits are independent of \(W\), although both the
numerator and denominator depend strongly on it.  Numerically, the weighted
nonbacktracking matrix has
\[
   \rho_{B,W}\approx 2.164666714,
\]
whereas the unweighted matrix has spectral radius \(2\).

The following values are exact integer weighted sums.  They were computed by the
weighted recurrence in \ref{app:dynamic-count}; the denominator was also
checked independently as \(\Tr(B_W^{\ell})\).
\begin{center}
\resizebox{0.82\textwidth}{!}{%
\(
\begin{array}{c|r|r|c}
\ell & N_{W,\ell}(K_4) & N_{W,\ell}^{6\text{-}\NZ}(K_4)
& N_{W,\ell}^{6\text{-}\NZ}(K_4)/N_{W,\ell}(K_4)\\ \hline
30 & 11527359956 & 3864928080 & 0.335283\\
31 & 24929521534 & 5594388712 & 0.224408\\
32 & 54043416200 & 18090312832 & 0.334737\\
33 & 116956641172 & 26150197524 & 0.223589\\
34 & 252980643848 & 84563637744 & 0.334269\\
35 & 547935284974 & 122344002620 & 0.223282\\
36 & 1186261793216 & 396554107752 & 0.334289\\
37 & 2566754539840 & 572117970376 & 0.222895\\
38 & 5557040394940 & 1856715479536 & 0.334119\\
39 & 12031254839686 & 2682390330360 & 0.222952\\
40 & 26038588031656 & 8698088979040 & 0.334046
\end{array}
\)}
\end{center}
The even and odd rows approach \(1/3\) and \(2/9\), respectively.  Hence the
canonical parity modulation persists under a genuinely directed perturbation of
the walk weights, while its two leading densities remain the graph invariants
specified by \(C_{K_4}(6)\) and \(D_{K_4}(6)\).
\end{example}

\begin{example}[Period modulation on \(\Theta_{2,2,4}\)]
\label{ex:bipartite-period-modulation}
Let \(G=\Theta_{2,2,4}\) be the theta graph formed by three internally disjoint
\(v\)-\(w\) paths \(\Delta_0,\Delta_1,\Delta_2\) of lengths \(2,2,4\),
respectively.  This is the graph \(G_3\) in the spectral-radius example
of~\cite[Subsection~4.1]{LuoRoyHSGT}.

\begin{figure}[H]
\centering
\begin{tikzpicture}[scale=0.85]
\coordinate (v) at (0,0);
\coordinate (w) at (6,0);
\coordinate (a) at (3,1.5);
\coordinate (b) at (3,0);
\coordinate (c1) at (1.5,-1.3);
\coordinate (c2) at (3,-1.7);
\coordinate (c3) at (4.5,-1.3);
\draw (v)--(a)--(w);
\draw (v)--(b)--(w);
\draw (v)--(c1)--(c2)--(c3)--(w);
\foreach \x in {v,w,a,b,c1,c2,c3} \fill (\x) circle (2pt);
\node[left] at (v) {\(v\)};
\node[right] at (w) {\(w\)};
\node at (3,1.05) {\(\Delta_0\)};
\node at (3,0.35) {\(\Delta_1\)};
\node at (3,-2.05) {\(\Delta_2\)};
\end{tikzpicture}
\caption{The bipartite theta graph \(\Theta_{2,2,4}\).  The three
\(v\)-\(w\) paths have lengths \(2,2,4\).}
\label{fig:theta-224}
\end{figure}
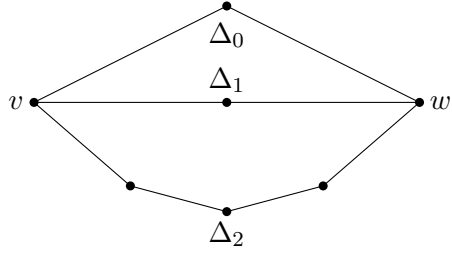

The graph is bipartite, has genus \(g=2\), and has nonbacktracking period
\[
   \nu_G=\gcd\{4,6\}=2.
\]
Orient every path from \(v\) to \(w\).  A mod-\(t\) flow is constant on each
path, and its three path values satisfy
\[
   x_0+x_1+x_2=0.
\]
Choosing nonzero \(x_0,x_1\) leaves \(x_2\ne0\) unless
\(x_1=-x_0\), and therefore
\[
   \CG(t)=(t-1)(t-2).
\]

The period character exists and, on integral flows in these path coordinates, is
\[
   \eta(x_0,x_1,x_2)=(-1)^{x_2},
\]
as computed in~\cite[Subsection~4.1]{LuoRoyHSGT}.  For positive even \(t\), its
mod-\(t\) descent \(\eta_t\) is therefore given by the same formula.  Fixing a
nonzero value of \(x_2\) leaves exactly \(t-2\) choices of nonzero
\((x_0,x_1)\), and Proposition~\ref{prop:period-coefficient} gives
\[
   \Deta(t)
      =(t-2)\sum_{x_2=1}^{t-1}(-1)^{x_2}
      =2-t.
\]
Since both sides are polynomials and this identity holds at every positive even
integer, it is a polynomial identity.  At \(t=4\), the nowhere-zero flows are the
permutations of \((1,1,2)\) and \((2,3,3)\).  Two have \(\eta=1\) and four have
\(\eta=-1\), so
\[
   \CG(4)=6,
   \qquad
   \Deta(4)=-2.
\]

The normalized period parity is also visible directly.  Since every internal
vertex of each \(\Delta_i\) has degree two, a circuit decomposes into complete
traversals of the three \(v\)-\(w\) paths.  Let \(q(P)\) be the total number of
these path traversals and let \(r(P)\) be the number that traverse \(\Delta_2\).
Because \(P\) is closed, \(q(P)\) is even, and
\[
   |P|=2\bigl(q(P)-r(P)\bigr)+4r(P)=2\bigl(q(P)+r(P)\bigr).
\]
Consequently, for every circuit of length \(\ell\),
\[
   \eta(P)=(-1)^{r(P)}=(-1)^{\ell/2}.
\]
Choose a \(\{\pm1\}\)-valued gain \(\gamma_\eta\) representing \(\eta\).  Then
\[
   \Tr(B_{\gamma_\eta}^\ell)=(-1)^{\ell/2}N_\ell(G).
\]
This is the specialization of the spectral rotation in
Proposition~\ref{prop:spectral-input}(c).  Theorem~\ref{thm:main}(b) therefore gives,
along even lengths,
\[
   N_{\ell}^{4\text{-}\NZ}(G)
      =
      \frac{6-2(-1)^{\ell/2}}{16}N_\ell(G)+O(\sigma^\ell)
\]
for some \(\sigma<\rho_B\).
It follows that
\[
   \lim_{\substack{\ell\to\infty\\ \ell\equiv0\ (\mathrm{mod}\ 4)}}
      \frac{N_{\ell}^{4\text{-}\NZ}(G)}{N_\ell(G)}
      =\frac14,
   \qquad
   \lim_{\substack{\ell\to\infty\\ \ell\equiv2\ (\mathrm{mod}\ 4)}}
      \frac{N_{\ell}^{4\text{-}\NZ}(G)}{N_\ell(G)}
      =\frac12.
\]
Since \(G\) is bipartite, \(D_G(4)=C_G(4)=6\), but the edge correction is governed
by \(D_{G,\eta}(4)=-2\), not by \(D_G(4)\).  The uniform density
\(C_G(4)/4^2=3/8\) is exactly the average of the two period-parity limits, while
\(D_{G,\eta}(4)\) measures their difference.  For contrast, the bipartite theta
graph \(\Theta_{1,3,5}\), denoted \(G_2\) in~\cite[Subsection~4.1]{LuoRoyHSGT},
also has period two but has no period character.  Its nonbacktracking formula has
only the ordinary flow-polynomial term.
\end{example}

\begin{example}[Evaluation-specific zeros]\label{ex:evaluation-zeros}
Evaluating the subset expansions defining \(C_G\) and \(D_G\) gives, for the
Petersen graph,
\[
   \CG(t)
      =
      (t-1)(t-2)(t-3)(t-4)(t^2-5t+10)
\]
and
\[
   \DG(t)=-5(t-2)(t-4)(t-6).
\]
Thus
\[
   \CG(2)=\CG(4)=0,
   \qquad
   \DG(2)=\DG(4)=0.
\]
At \(t=6\), however,
\[
   \CG(6)=1920,
   \qquad
   \DG(6)=0.
\]
Thus the Petersen graph has many nowhere-zero \(6\)-flows, but no leading parity
bias at \(t=6\).

For a second comparison, let \(H\) be obtained from \(K_5\) by deleting two
adjacent edges.  The same calculation gives
\[
   C_H(t)=(t-2)(t-1)(t^2-4t+5),
   \qquad
   D_H(t)=(t-2)(t-5).
\]
At \(t=5\),
\[
   C_H(5)=120,
   \qquad
   D_H(5)=0.
\]
At \(t=4\), the same graph has
\[
   C_H(4)=30,
   \qquad
   D_H(4)=-2,
\]
and hence
\[
   r_{\mathrm{even}}=\frac{C_H(4)+D_H(4)}{4^4}=\frac{7}{64},
   \qquad
   r_{\mathrm{odd}}=\frac{C_H(4)-D_H(4)}{4^4}=\frac{1}{8}.
\]
Thus \(\DG(t)\) can vanish at a fixed integer even when the ordinary flow
polynomial is positive there, while recording a nonzero parity bias at a neighboring
evaluation.  In both graphs, the displayed zero is an evaluation-specific
cancellation in a nonzero polynomial, not an obstruction to nowhere-zero flows.
\end{example}

The canonical examples are summarized as follows:
\[
\begin{array}{c|c|c|c|c|c}
G & t & \CG(t) & \DG(t) & r_{\mathrm{even}} & r_{\mathrm{odd}}\\ \hline
K_4 & 4 & 6 & 6 & 3/16 & 0\\
K_4 & 5 & 24 & 9 & 24/125 & 24/125\\
K_4 & 6 & 60 & 12 & 1/3 & 2/9\\
\mathrm{Petersen} & 6 & 1920 & 0 & 10/243 & 10/243\\
H & 4 & 30 & -2 & 7/64 & 1/8
\end{array}
\]
Here \(H=K_5-\{2\ \text{adjacent edges}\}\), as in
Example~\ref{ex:evaluation-zeros}.  The period example has two admissible length
classes rather than an even--odd split and is summarized separately:
\[
\begin{array}{c|c|c|c|c|c}
G & t & \CG(t) & \Deta(t)
  & r_{\ell\equiv0\ (\mathrm{mod}\ 4)}
  & r_{\ell\equiv2\ (\mathrm{mod}\ 4)}\\ \hline
\Theta_{2,2,4} & 4 & 6 & -2 & 1/4 & 1/2
\end{array}
\]
For the odd evaluation \(t=5\), the even- and odd-length canonical densities are
both equal to \(C_{K_4}(5)/5^3\).  Thus
\(K_4\) gives complete parity concentration at \(t=4\), a single limiting density
at \(t=5\), and a nonzero two-level parity oscillation at \(t=6\).  The Petersen
graph gives no leading parity bias at \(t=6\), while \(H\) gives a nonzero bias
with both parity classes present.

For the nonbipartite graphs above, the degrees of \(D_G(t)\) are \(1\) for
\(K_4\), \(3\) for the Petersen graph, and \(2\) for \(H\).  In each case the
degree is strictly smaller than \(g=\beta(G)\), consistently with the bound in
Proposition~\ref{prop:canonical-coefficient}.  These graphs are also bridgeless,
and their displayed polynomials exhibit the strict coefficient alternation asserted by
Proposition~\ref{prop:semimatroid-specialization}.

\appendix
\renewcommand{\thetheorem}{\Alph{section}.\arabic{theorem}}

\section{Dynamic counting of circuits}\label{app:dynamic-count}

This appendix records the exact enumeration used in
Example~\ref{ex:k4-counts}.  The point is not to introduce a new theoretical
ingredient, but to make clear that the numbers in the example are obtained by
direct counting of rooted, oriented circuits, with their
mod-\(t\) abelianization tracked explicitly.

Fix an orientation \(\mathfrak o\) of
\[
   E(G)=\{e_1,\ldots,e_m\}.
\]
For an oriented edge \(\mathbf a\in\mathbf E(G)\), define
\[
   \mathbf s(\mathbf a)\in(\bbZ/t\bbZ)^E
\]
by
\[
   \mathbf s(\mathbf a)_{e_i}
      =
      \begin{cases}
      1, & \mathbf a=\mathbf e_i,\\
      -1,& \mathbf a=\mathbf e_i^{-1},\\
      0, & \text{otherwise}.
      \end{cases}
\]
Thus \(\mathbf s(\mathbf a)\) is the contribution of one traversal of
\(\mathbf a\) to the mod-\(t\) abelianization.  Write
\[
   \mathbf a\leadsto\mathbf b
\]
when \(\mathbf a(1)=\mathbf b(0)\) and \(\mathbf b\ne \mathbf a^{-1}\).  This is the
ordinary nonbacktracking transition relation.

For \(h\ge1\), let
\[
   A_h(\mathbf a,\mathbf b,r)
\]
be the number of nonbacktracking walks of length \(h\) whose first oriented edge is
\(\mathbf a\), whose last oriented edge is \(\mathbf b\), and whose accumulated
residue vector is \(r\in(\bbZ/t\bbZ)^E\).  The initial values and recurrence are
\[
   A_1(\mathbf a,\mathbf b,r)
      =
      \mathbf 1_{\{\mathbf b=\mathbf a\}}\,
      \mathbf 1_{\{r=\mathbf s(\mathbf a)\}},
\]
and, for \(h\ge1\),
\[
   A_{h+1}(\mathbf a,\mathbf c,r)
      =
      \sum_{\mathbf b:\,\mathbf b\leadsto\mathbf c}
         A_h\bigl(\mathbf a,\mathbf b,r-\mathbf s(\mathbf c)\bigr).
\]
Let
\[
   \mathcal R_t^{\NZ}
      =
      \{r\in(\bbZ/t\bbZ)^E:\ r_e\ne0\ \text{for every }e\in E(G)\}.
\]
Then the exact count used in Example~\ref{ex:k4-counts} is
\[
   N_{\ell}^{t\text{-}\NZ}(G)
      =
      \sum_{\substack{\mathbf a,\mathbf b\in\mathbf E(G)\\
                      \mathbf b\leadsto\mathbf a}}
      \sum_{r\in\mathcal R_t^{\NZ}}
         A_\ell(\mathbf a,\mathbf b,r).
\]
For a positive directed weight \(W\), define instead
\[
   A_{W,1}(\mathbf a,\mathbf b,r)
      =
      W(\mathbf a)\,
      \mathbf 1_{\{\mathbf b=\mathbf a\}}\,
      \mathbf 1_{\{r=\mathbf s(\mathbf a)\}},
\]
and
\[
   A_{W,h+1}(\mathbf a,\mathbf c,r)
      =
      W(\mathbf c)
      \sum_{\mathbf b:\,\mathbf b\leadsto\mathbf c}
         A_{W,h}\bigl(\mathbf a,\mathbf b,r-\mathbf s(\mathbf c)\bigr).
\]
Every extension is thereby multiplied by the weight of its newly appended edge,
so induction gives
\[
   N_{W,\ell}^{t\text{-}\NZ}(G)
      =
      \sum_{\substack{\mathbf a,\mathbf b\in\mathbf E(G)\\
                      \mathbf b\leadsto\mathbf a}}
      \sum_{r\in\mathcal R_t^{\NZ}}
         A_{W,\ell}(\mathbf a,\mathbf b,r).
\]
Removing the restriction \(r\in\mathcal R_t^{\NZ}\) gives \(N_{W,\ell}(G)\).
This weighted form of the recurrence produces the exact values in
Example~\ref{ex:weighted-k4}; setting \(W\equiv1\) recovers the preceding
unweighted recurrence.

Indeed, by induction on \(h\), \(A_h(\mathbf a,\mathbf b,r)\) counts precisely the
nonbacktracking paths with the prescribed first edge, last edge, and residue.  The
condition \(\mathbf b\leadsto\mathbf a\) imposes the final cyclic nonbacktracking
condition, so the sum above counts the same rooted, oriented circuits as
\(\Tr(B_0^\ell)\), retaining only those whose mod-\(t\) abelianization is nowhere
zero on every edge.  If the restriction \(r\in\mathcal R_t^{\NZ}\) is removed, the
same formula gives the total nonbacktracking count \(N_\ell(G)\).

In the accompanying notebook, the residue vector \(r\) is encoded by a base-\(t\)
integer and the updates \(r\mapsto r+\mathbf s(\mathbf c)\) are precomputed as
permutations of the residue space.  This is only an implementation of the recurrence
above; all counts are exact integer counts, with no use of random sampling or of the
asymptotic formula from Theorem~\ref{thm:main}.

\section{A Perron--Frobenius proof for postman circuits}
\label{app:postman-pf}

We give a second proof of Corollary~\ref{cor:postman-asymptotic}.  Unlike the proof
in Subsection~\ref{sec:weighted-postman}, it does not pass through integral
abelianization.  It begins with an
exact weighted extension of the covering-cycle formula of da Costa and
Policarpo~\cite[Theorem~2.1]{DaCostaPolicarpo2018}, and consequently applies to
graphs with bridges as long as no tree is attached to the cyclic core.

For \(F\subseteq E(G)\), let \(B_{F,W}\) be the weighted nonbacktracking edge
adjacency matrix of the spanning subgraph \(G[F]\), with \(W\) restricted to the
oriented edges above \(F\).  Isolated vertices play no role in this matrix.

\begin{proposition}[Weighted covering-cycle formula]
\label{prop:postman-inclusion-exclusion}
For every finite graph \(G\), every positive directed weight \(W\), and every
\(\ell\ge1\),
\begin{equation}
   \pc_{W,\ell}(G)
      =
      \sum_{F\subseteq E(G)}
         (-1)^{|E(G)|-|F|}\Tr(B_{F,W}^{\ell}).
   \label{eq:postman-inclusion-exclusion}
\end{equation}
\end{proposition}

\begin{proof}
Expand the right-hand side as a weighted sum over rooted, oriented
circuits.  A circuit \(P\) with unoriented edge support
\(S\subseteq E(G)\) occurs in \(\Tr(B_{F,W}^{\ell})\) exactly when \(S\subseteq F\).
Its total coefficient in \eqref{eq:postman-inclusion-exclusion} is therefore
\[
   \sum_{S\subseteq F\subseteq E(G)}(-1)^{|E(G)|-|F|}
      =
      \begin{cases}
      1,&S=E(G),\\
      0,&S\ne E(G).
      \end{cases}
\]
Thus precisely the postman circuits remain, each with its original weight \(W(P)\).
\end{proof}

Formula~\eqref{eq:postman-inclusion-exclusion} is the rooted weighted form of
da Costa and Policarpo's unweighted formula~\cite[Theorem~2.1]{DaCostaPolicarpo2018}.
A directed-support transform based on the same inclusion-exclusion principle was
later obtained independently in
\cite[Appendix~A]{BandExnerGoelStrauss2026} for periodic-orbit and Eulerian-circuit
counting.

The transfer-operator viewpoint behind the spectral argument below has a broader
precedent in the work of Adachi and Sunada~\cite{AdachiSunada1987}, who developed
twisted Perron--Frobenius theory for topological graphs and the associated
\(L\)-functions.  In the present finite-dimensional setting, the strict comparison
needed here follows directly from the equality case of Wielandt's matrix theorem.

\begin{proposition}[Direct spectral estimate]
\label{prop:postman-pf}
Let \(G\) be connected, have genus \(g\ge2\), and have no vertex of degree one.
For every positive directed weight \(W\), there is a number
\(\sigma<\rho_{B,W}\) such that
\[
   \pc_{W,\ell}(G)=N_{W,\ell}(G)+O(\sigma^\ell).
\]
Consequently,
\[
   \pc_{W,\ell}(G)
      =
      \begin{cases}
      \nu_G\rho_{B,W}^{\ell}+O(\sigma^\ell),&\nu_G\mid\ell,\\
      0,&\nu_G\nmid\ell.
      \end{cases}
\]
\end{proposition}

\begin{proof}
Write \(B=B_W\) and \(\rho=\rho_{B,W}\).  We first verify that the
Perron--Frobenius theorem applies to \(B\).  Since \(G\) is connected and has no
degree-one vertex, it is its own cyclic core.  The assumption \(g\ge2\) excludes a
cycle.  The elementary nonbacktracking-connectivity lemma
\cite[Lemma~5.9]{LuoRoyHSGT}, whose relevant assertion goes back to Ortner and
Woess~\cite[Lemma~3.1]{OrtnerWoess2007}, therefore shows that, for every pair of
oriented edges \(\mathbf a,\mathbf b\in\mathbf E(G)\), there is a nonbacktracking
edge-walk from \(\mathbf a\) to \(\mathbf b\).  The support digraph of \(B\) is
consequently strongly connected.  Since all directed weights are positive, they do
not alter this support, and hence \(B\) is an irreducible nonnegative matrix.

We next prove the strict spectral-radius estimate needed for every proper edge
set.  Fix \(F\subsetneq E(G)\), and extend \(B_{F,W}\) to a matrix
\(\widehat B_{F,W}\) indexed by all of \(\mathbf E(G)\) by putting
\[
   (\widehat B_{F,W})_{\mathbf a,\mathbf b}
      =
      \begin{cases}
      B_{\mathbf a,\mathbf b},
         &\mathbf a,\mathbf b\in\mathbf E(G[F]),\\
      0,&\text{otherwise}.
      \end{cases}
\]
After a simultaneous permutation of rows and columns,
\(\widehat B_{F,W}=B_{F,W}\oplus0\), and therefore
\[
   \rho(\widehat B_{F,W})=\rho(B_{F,W}).
\]
Moreover,
\[
   0\le \widehat B_{F,W}\le B
\]
entrywise, and the two matrices are unequal because \(F\ne E(G)\) and the
irreducible matrix \(B\) has no zero row.

We now apply Wielandt's equality theorem~\cite{Wielandt1950}; see also
\cite[Theorem~8.4.5]{HornJohnson2013}.  It states that if \(A\) is nonnegative and
irreducible, \(|M|\le A\), and \(\rho(M)=\rho(A)\), then
\[
   M=\ee^{\mathrm i\vartheta}DAD^{-1}
\]
for some real \(\vartheta\) and some diagonal unitary matrix \(D\).  In particular,
equality of the spectral radii forces \(|M|=A\) entrywise.  Apply this with
\(A=B\) and \(M=\widehat B_{F,W}\).  If their spectral radii were equal, then,
because both matrices are nonnegative, we would obtain
\[
   \widehat B_{F,W}=|\widehat B_{F,W}|=B,
\]
contrary to the proper-compression observation above.  Consequently,
\[
   \rho(B_{F,W})<\rho
   \qquad(F\subsetneq E(G)).
\]

Because there are finitely many edge subsets, the number
\[
   \tau=
   \max\bigl(\{0\}\cup
      \{\rho(B_{F,W}):\varnothing\ne F\subsetneq E(G)\}\bigr)
\]
satisfies \(\tau<\rho\).  If
\(\lambda_{F,1},\ldots,\lambda_{F,2|F|}\) are the eigenvalues of \(B_{F,W}\),
listed with algebraic multiplicity, then
\[
   \left|\Tr(B_{F,W}^{\ell})\right|
      =
      \left|\sum_{j=1}^{2|F|}\lambda_{F,j}^{\ell}\right|
      \le 2|F|\rho(B_{F,W})^{\ell}
      \le 2|E(G)|\tau^{\ell}.
\]
The empty subgraph contributes zero.  Summing over the finitely many proper
subsets and separating the term \(F=E(G)\) in
\eqref{eq:postman-inclusion-exclusion} gives
\[
   \pc_{W,\ell}(G)
      =\Tr(B^{\ell})+O(\tau^{\ell})
      =N_{W,\ell}(G)+O(\tau^{\ell}).
\]

To evaluate the leading term, put
\(\zeta=\exp(2\pi\mathrm i/\nu_G)\).  The peripheral-spectrum part of the
Perron--Frobenius theorem states that, for an irreducible matrix of period
\(\nu_G\), the eigenvalues of modulus \(\rho\) are precisely
\[
   \rho,\ \rho\zeta,\ \ldots,\ \rho\zeta^{\nu_G-1},
\]
each with algebraic multiplicity one.  All remaining eigenvalues have modulus at
most some \(\sigma_0<\rho\).  Therefore
\[
   \Tr(B^{\ell})
      =
      \rho^{\ell}\sum_{j=0}^{\nu_G-1}\zeta^{j\ell}
         +O(\sigma_0^{\ell}).
\]
When \(\nu_G\mid\ell\), the sum is \(\nu_G\), giving
\[
   N_{W,\ell}(G)
      =\nu_G\rho^{\ell}+O(\sigma_0^{\ell}).
\]
When \(\nu_G\nmid\ell\), no closed edge-walk of length \(\ell\) occurs in an
irreducible shift of period \(\nu_G\); hence
\(N_{W,\ell}(G)=\pc_{W,\ell}(G)=0\) exactly.  Taking
\(\sigma=\max\{\tau,\sigma_0\}<\rho\) completes the proof.
\end{proof}

The direct proof shows that the bridgelessness assumption in
Corollary~\ref{cor:postman-asymptotic} is not intrinsic to postman circuits.  It is
needed there only because integral nowhere-zero abelianization is impossible on a
bridge.  A bridge joining cyclic subgraphs may occur in a postman circuit, once in
each direction, and is covered by Proposition~\ref{prop:postman-pf}.

\end{document}